\documentclass[12pt]{amsart}

\usepackage{amsfonts,amsmath,amssymb,amsthm,mathrsfs,url,xcolor,latexsym,rawfonts,graphicx}
\usepackage{hyperref}
\usepackage{enumerate} 
\usepackage[]{caption2}
\usepackage{accents}
\usepackage{amssymb,amsmath}
\usepackage{amsbsy}
\usepackage{cite}
\usepackage{float}
\usepackage{graphicx}

\theoremstyle{plain}

   \newcommand{\beqn}{\begin{eqnarray}}
   \newcommand{\eeqn}{\end{eqnarray}}
   \newcommand{\beqs}{\begin{eqnarray*}}
   \newcommand{\eeqs}{\end{eqnarray*}}
   \newcommand{\ban}{\begin{eqnarray*}}
   \newcommand{\nan}{\end{eqnarray*}}
   \newcommand{\beq}{\begin{equation}}
   \newcommand{\eeq}{\end{equation}}

\renewcommand{\det}{\mbox{det}}

\numberwithin{equation}{section}

  \numberwithin{equation}{section}
  \numberwithin{figure}{section}

\begin{document}

\title[The Gaussian Minkowski problem]
{\textbf{Existence of non-symmetric solutions to the Gaussian Minkowski problem}}

\author[Y. Feng]
{Yibin Feng}
\address{Yibin Feng, School of Mathematical Sciences,
University of Science and Technology of China,
Hefei, 230026, P.R. China.}
\email{fybt1894@ustc.edu.cn}

\author[W. Liu]
{Weiru Liu}
\address
	{Weiru Liu, School of Mathematical Sciences,
University of Science and Technology of China,
Hefei, 230026, P.R. China.}
\email{lwr19997@ustc.edu.cn}

\author[L. Xu]
{Lei Xu}
\address
{Lei Xu, School of Mathematical Sciences,
University of Science and Technology of China,
Hefei, 230026, P.R. China.}
\email{xlsx@mail.ustc.edu.cn}



\thanks{This work was supported by the China Postdoctoral Science Foundation (Grant No.2020TQ0315). 
 }


\date{\today}

\dedicatory{}

\begin{abstract}  
Existence of symmetric solutions to the Gaussian Minkowski problem was established by Huang, Xi and Zhao. In this paper, we show the existence of non-symmetric solutions to this problem by studying the related Monge-Amp\`{e}re type equation on the sphere.
\end{abstract}

\maketitle

\baselineskip=16.4pt
\parskip=3pt

\section{\bf Introduction}
 Let $\mathbb{R}^n$ be the $n$-dimensional Euclidean space. The unit sphere in $\mathbb{R}^n$ is denoted by $\mathbb{S}^{n-1}$. A convex body in $\mathbb{R}^n$ is a compact convex set with non-empty interior. Denote by $\mathcal{K}_o^n$ the set of all convex bodies that contain the origin in their interiors in $\mathbb{R}^n$ and by $\mathcal{K}_e^n$ the set of all origin-symmetric convex bodies in $\mathbb{R}^n$.
 
 The Gaussian volume, denoted by $\gamma_n$, of $K\in\mathcal{K}_o^n$ is given by
 \begin{eqnarray}\label{1.0}
 	\gamma_n(K)=\frac{1}{\left(\sqrt{2\pi}\right)^n}\int_Ke^{-\frac{|y|^2}{2}}dy.
 \end{eqnarray}
For $K, L\in\mathcal{K}_o^n$, the Gaussian surface area measure $S_{\gamma_n, K}$, discovered in \cite{1a}, can be defined by the variational formula
 \begin{eqnarray*}
 	\lim_{t\rightarrow 0}\dfrac{\gamma_n(K+tL)-\gamma_n(K)}{t}=\int_{\mathbb{S}^{n-1}}h_LdS_{\gamma_n, K},
 \end{eqnarray*}
where $K+tL=\{y+tz:y\in K ~\mbox{and}~ z\in L\}$ is the Minkowski combination of $K$ and $L$, and $h_L:\mathbb{S}^{n-1}\rightarrow \mathbb{R}$ is the support function of $L$. If $K\in\mathcal{K}_o^n$ has a $C^2$ boundary with positive curvature everywhere, then $S_{\gamma_n, K}$ is absolutely continuous with respect to the spherical Lebesgue measure, and its density is equal to
 \begin{eqnarray}\label{1.1}
\frac{1}{(\sqrt{2\pi})^n}e^{-\frac{\left|\nabla h_K\right|^2+h^2_K }{2}}\det\left(\nabla^2h_K+h_KI\right)~\mbox{on}~\mathbb{S}^{n-1}, 
\end{eqnarray}
where $h_K:\mathbb{S}^{n-1}\rightarrow \mathbb{R}$ is the support function of $K$, $\nabla$ and $\nabla^2$ are gradient and Hessian operators with respect to an orthonormal frame on $\mathbb{S}^{n-1}$, and $I$ is the identity matrix.

The Gaussian Minkowski problem is the counterpart of the classical Minkowski problem characterizing surface area measure. The classical Minkowski problem asks for necessary and sufficient conditions so that a given measure on the unit sphere is the surface area measure of a convex body. Its solution amounts to solving a fully nonlinear partial differential equaton when a given measure has a density function. The study of the problem has a long history and influence many fields of mathematics including convex geometric analysis, partial differential equaton, functional analysis, and so on. The Gaussian Minkowski problem posed in \cite{1a} is the problem characterizing Gaussian surface area measure. Namely,
 
 \noindent{\bf The Gaussian Minkowski problem.} Given a finite Borel measure $\mu$ on the unit sphere $\mathbb{S}^{n-1}$, what are the necessary and sufficient conditions on $\mu$ so that there exists a convex body $K\in\mathcal{K}_o^n$ such that
\begin{eqnarray}\label{1.2}
S_{\gamma_n, K}=\mu?
\end{eqnarray} 
It follows from (\ref{1.1}) that when the measure $\mu$ has a density function $\frac{1}{f}$, the Gaussian Minkowski problem is equivalent to the solvability of following Monge-Amp\`{e}re type equation on the sphere:
 \begin{eqnarray}\label{1.3}
	\frac{1}{(\sqrt{2\pi})^n}e^{-\frac{\left|\nabla h_K\right|^2+h^2_K }{2}}\det\left(\nabla^2h_K+h_KI\right)=\frac{1}{f}.
\end{eqnarray}
 
Huang, Xi and Zhao \cite{1a} studied the even Gaussian Minkowski problem and first obtained the following existence result.
 
{\it \noindent{\bf Theorem A \cite{1a}.}~Suppose $\mu$ is a finite even Borel measure not concentrated in any closed hemisphere. Then there exist an origin-symmetric convex body $K\in \mathcal{K}_e^n$ and $\tau>0$ such that 
\begin{eqnarray*}
S_{\gamma_n, K}=\tau\mu.
\end{eqnarray*}} 
Theorem A was proved by the variational method in \cite{1a}. Since the Gaussian surface area measure does not have any homogeneity, the constant $\tau$ can not be removed in general. Likewise, solutions to the Orlicz Minkowski problem usually also contain a constant that is not equal to $1$. See for example \cite{2a, 3a, 4a, 5a, 6a, 7a, 8a, 9a, 10a, 11a}. These solutions with the constant factors are often referred to as the normalization solutions for this type of Minkowski problems involving non-homogeneous geometric and probabilistic measures.
 
The first aim of this paper is to give the following normalization solution to the Gaussian Minkowski problem for general measures. That is, we get rid of the even assumption on $\mu$ in Theorem A. Our result is stated below.
 
 {\it \noindent{\bf Theorem 1.1.}~Suppose $\mu$ is a finite Borel measure not concentrated in any closed hemisphere. Then there exist a convex body $K\in \mathcal{K}_o^n$ and a constant $\tau>0$ such that 
 	\begin{eqnarray*}
 		S_{\gamma_n, K}=\tau\mu.
 \end{eqnarray*}} 
Theorem 1.1 was also considered in \cite{12a}, where a variational argument was used by establishing the variational formula and the optimizaiton problem associated to the $L_p$-Gaussian Minkowski problem. Our approach is based on the study of a suitably parabolic flow and the use of approximation argument. The main difficulties of proving Theorem 1.1 is to obtain uniformly positive lower and upper bounds for support function and principal curvatures along the flow. By choosing proper auxiliary functions and initial value, the bounds of support function and principal curvatures are obtained after careful analysis and computations.
 
As the solutions to the Orlicz Minkowski problem, it is very difficult to remove the constant factor $\tau$ in Theorem A because of the lack of homogeneity for the Gaussian surface area measure. In \cite{1a}, Huang, Xi and Zhao used a degree-theoretic approach combining the approximation argument to overcome this difficulty, and obtained a non-normalized solution (meaning that there is no constant $\tau$) to the even Gaussian Minkowski problem. Namely, they established the following result.
 
{\it \noindent{\bf Theorem B \cite{1a}.}~Let $\mu$ be an  even measure on $\mathbb{S}^{n-1}$ that is not concentrated in any subspace and $|\mu|<\frac{1}{\sqrt{2\pi}}$. Then there exists a unique origin-symmetric convex body $K\in \mathcal{K}_e^n$ with $\gamma_n(K)>\frac{1}{2}$ such that 
\begin{eqnarray*}
		S_{\gamma_n, K}=\mu.
\end{eqnarray*}}  
The second aim of this paper is to give the following existence result for Theorem B without symmetry hypothesis on the given measure $\mu$. 

{\it \noindent{\bf Theorem 1.2.}~Let $\mu$ be a non-zero  finite Borel measure on $\mathbb{S}^{n-1}$ that is not concentrated in any closed hemisphere. Then for $|\mu|<\frac{1}{\sqrt{2\pi}}$, there exists a  convex body $K\in \mathcal{K}_o^n$ with $\gamma_n(K)>\frac{1}{2}$ such that 
	\begin{eqnarray*}
		S_{\gamma_n, K}=\mu.
\end{eqnarray*}}   
Building upon the above existence result, we first convert the original Monge-Amp\`{e}re type equation (\ref{1.3}) on the unit sphere $\mathbb{S}^{n-1}$ into a Euclidean Monge-Amp\`{e}re type equation on $\mathbb{R}^{n-1}$, then study the regularity of the solutions to the Gaussian Minkowski problem based on Caffarelli's results \cite{13a, 14a}, further get the existence of smooth solutions to the Gaussian Minkowski problem following the degree-theoretic argument in \cite{1a, 12a}, and finally prove Theorem 1.2 by an approximation argument. 
 
Since the seminal work of Huang, Xi and Zhao \cite{1a}, the Gaussian Minkowski problem has quickly become the topic of many influential works; see for example \cite{12a, 15a, 16a}.
 
 The rest of the paper is organized as follows. In Section 2, we give some preliminaries that will be used in the proofs of our results. In Section 3, the geometric flow and its associated functional will be introduced, the priori estimates will be demonstrated to obtain the existence of smooth solutions to the Gaussian Minkowski problem, and the proof of Theorem 1.1 will be completed which implies the existence of weak solutions to the Gaussian Minkowski problem. In Section 4, we will be devoted to proving Theorem 1.2 based on the famous regularity results by Caffarelli \cite{13a, 14a} and the degree-theoretic method.

 \section{\bf Preliminaries}

\subsection{\bf Convex bodies and Gaussian surface area measure}\hfill \break
\indent
For quick later consultation we collect some notation and basic facts about convex bodies and Gaussian surface area measure. Good references are the books by Gardner \cite {17a} and Schneider\cite{18a} as well as the papers\cite{19a} and \cite{1a}.
 
The standard inner product of the vectors $y, z\in \mathbb{R}^n$ is denoted by $\langle y, z\rangle$. We write $|y|=\sqrt{\langle y, y\rangle}$. The boundary of a convex body $K$ is written as $\partial K$. We use  $\mathcal{H}^{k}$ to denote the $k$-dimensional Hausdorff measure. The set of continuous functions on the unit sphere $\mathbb{S}^{n-1}$ will be denoted by $C(\mathbb{S}^{n-1})$ and will always be viewed as equipped with the max-norm metric
\begin{eqnarray*}
\left\|f-g\right\|_\infty=\max_{x\in \mathbb{S}^{n-1}}\left| f(x)-g(x)\right|, 
\end{eqnarray*}
for $f, g\in C(\mathbb{S}^{n-1})$. The subspace of positive continuous functions from $C(\mathbb{S}^{n-1})$ will be denoted by $C^+(\mathbb{S}^{n-1})$.

We will need the trivial fact, cf.\cite{7a}, that for $f, f_i\in C(\mathbb{S}^{n-1})$ such that $f_i\rightarrow f$ uniformly on $\mathbb{S}^{n-1}$ and finite measures $\mu, \mu_i$ such that $\mu_i\rightarrow \mu$ weakly on $\mathbb{S}^{n-1}$, there is
 \begin{eqnarray}\label{2.1}
 	\lim_{i\rightarrow \infty}\int_{\mathbb{S}^{n-1}}f_i(x)d\mu_i(x)=\int_{\mathbb{S}^{n-1}}f(x)d\mu(x).
 \end{eqnarray}

For each $f\in C^+(\mathbb{S}^{n-1})$, the Wulff shape $\left[ f\right] $ generated by $f$ is the convex body defined by
 \begin{eqnarray*}
\left[ f\right]=\left\lbrace y\in \mathbb{R}^{n}: \langle y, x\rangle\leq f(x), ~\mbox{for all}~x\in \mathbb{S}^{n-1}\right\rbrace. 
\end{eqnarray*}
It is apparent that $h_{\left[ f\right]} \leq f $ and $\left[ h_K\right]=K$ for each $K\in \mathcal{K}_o^n$.

For a convex body $K\in \mathcal{K}_o^n$, the support function $h=h_K: \mathbb{R}^n\rightarrow \mathbb{R}$ and the radial function $\rho=\rho_K: \mathbb{R}^n\setminus\{o\}\rightarrow \mathbb{R}$ are respectively defined by
 \begin{eqnarray*}
	h(y)=\max\{\langle y, z\rangle: z\in K\},~~ \rho(z)=\max\{\lambda: \lambda z\in K\}.
\end{eqnarray*}
 The polar body of $K\in \mathcal{K}_o^n$ is defined by
\begin{eqnarray}\label{2.2}
	K^\ast=\{y\in\mathbb{R}^n: \langle y, z\rangle\leq 1~ \mbox {for all}~ z\in K\}.
\end{eqnarray}
From (\ref{2.2}), it easily follows that $K^\ast\in \mathcal{K}_o^n$ and $(K^\ast)^\ast=K$. In addition,
\begin{eqnarray}\label{2.3}
	\rho_K=1/h_{K^\ast}, ~~h_K=1/\rho_{K^\ast}.
\end{eqnarray}

It is easy to show that the support function of the line segment $\tilde{u}$ joining the points $o$ and $u\in \mathbb{S}^{n-1}$ is
\begin{eqnarray*}
	h_{\tilde{u}}(x)=\langle x, u\rangle_+=\max\{\langle x, u\rangle, 0\}, ~x\in \mathbb{S}^{n-1}.
\end{eqnarray*}

For $K\in \mathcal{K}_o^n$ with $\rho_K(u_{\min})=\min_{u\in\mathbb{S}^{n-1}}\rho_K(u)$ and $h_K(x_{\max})=\max_{x\in\mathbb{S}^{n-1}}h_K(x)$, the facts found in \cite{19a} that
\begin{eqnarray}\label{2.4}
	\max_{\mathbb{S}^{n-1}}h_K=\max_{\mathbb{S}^{n-1}}\rho_K,~~\min_{\mathbb{S}^{n-1}}h_K=\min_{\mathbb{S}^{n-1}}\rho_K, 
\end{eqnarray}
\begin{eqnarray}\label{2.5}
\rho_K(u)\langle u, u_{\min}\rangle\leq \rho_K(u_{\min}) ~~\mbox{for any}~~u\in \mathbb{S}^{n-1},
\end{eqnarray}
\begin{eqnarray}\label{2.5+}
	h_K(x)\geq \langle x, x_{\max}\rangle h_K(x_{\max}) ~~\mbox{for any}~~x\in \mathbb{S}^{n-1}
\end{eqnarray}
are of critical importance.

For a sequence of convex bodies $K_i \in\mathcal{K}_o^n$, we say that $K_i\rightarrow K_0\in \mathcal{K}_o^n$ with respect to the Hausdorff metric, provided when $i\rightarrow \infty$
\begin{eqnarray*}
	\|h_{K_i}-h_{K_0}\|_\infty=\max_{x\in \mathbb{S}^{n-1}}|h_{K_i}(x)-h_{K_0}(x)|\rightarrow 0,
\end{eqnarray*}
or equivalently
\begin{eqnarray*}
	\|\rho_{K_i}-\rho_{K_0}\|_\infty=\max_{u\in \mathbb{S}^{n-1}}|\rho_{K_i}(u)-\rho_{K_0}(u)|\rightarrow 0.
\end{eqnarray*}

For $y\in \partial K$, the spherical image of $y$ is defined as $\nu_K(y)=\{x\in\mathbb{S}^{n-1}: h_K(x)=\langle x, y\rangle\}$. For a Borel set $\eta \subset \mathbb{S}^{n-1}$, the reverse spherical image is $\nu_K^{-1}(\eta)=\{y\in \partial K: \nu_K(y)\cap \eta\neq \emptyset\}$. Let $\partial'K$ denote the subset of the boundary of $K$ where there is a unique outer unit normal vector. It is known that $\mathcal{H}^{n-1}(\partial K\backslash\partial'K)=0$. Then $\nu_K:\partial'K \rightarrow\mathbb{S}^{n-1}$ is a function that is usually called the spherical Gauss map.

The surface area measure $S_K$ of a convex body $K$ is a Borel measure on  $\mathbb{S}^{n-1}$ defined, for a Borel set $\omega\subset \mathbb{S}^{n-1}$, by
\begin{eqnarray*}
	S_K(\omega)=\mathcal{H}^{n-1}(\nu_K^{-1}(\omega))=\mathcal{H}^{n-1}(\{y\in \partial K: \nu_K(y)\cap \omega\neq \emptyset\}).
\end{eqnarray*}
 
It was shown in \cite{1a} that the Gaussian surface area measure of $K \in\mathcal{K}_o^n$ has the following integral representation
\begin{eqnarray}\label{2.6}
S_{\gamma_n, K}(\eta)=\frac{1}{\left( \sqrt{2\pi}\right) ^n}\int_{\nu_K^{-1}(\eta)}e^{-\frac{|y|^2}{2}}d\mathcal{H}^{n-1}(y)
\end{eqnarray}
for each Borel measure $\eta\subset \mathbb{S}^{n-1}$.

The Gaussian surface area measure is weakly convergent with respect to Hausdorff metric (see \cite{1a}): If $K_i \in\mathcal{K}_o^n$ such that $K_i\rightarrow K_0\in \mathcal{K}_o^n$ as $i\rightarrow\infty$, then $S_{\gamma_n, K_i}\rightarrow S_{\gamma_n, K_0}$ weakly.

The uniqueness of the solutions to the Gaussian Minkowski problem was established in \cite{1a}:

{\it\noindent{\bf Lemma 2.1\cite{1a}.}~Let $K, L\in \mathcal{K}_o^n$ such that $\gamma_n(K), \gamma_n(L)\geq \frac{1}{2}$. If $S_{\gamma_n, K}=S_{\gamma_n, L}$, then $K=L$. }

The Gaussian isoperimetric inequality implies the following result.

{\it\noindent{\bf Lemma 2.2\cite{1a}.}~If $K$ is a convex body in $\mathbb{R}^n$ such that $\gamma_n(K)=\frac{1}{2}$, then $\left|S_{\gamma_n, K} \right| \geq \frac{1}{\sqrt{2\pi}}$. }

\section{\bf Normalized solution}

In this section, Theorem 1.1 will be proved by the method of geometric flow. The geometric flow generated by the Gauss curvature was first studied by Firey \cite{20a} to describe the shape of a tumbling stone. Whereafter, the application of geometric flow in the relevant Minkowski problems has been extensively studied by many researches; see for example \cite{21a, 22a, 23a, 24a, 25a, 26a, 27a, 28a, 29a, 30a, 31a, 32a, 33a, 34a, 35a, 36a, 37a, 38a, 39a, 40a, 42a, 43a, 45a} and the references therein. Our geometric flow is inspired by the works of \cite{22a, 31a, 38a, 44a}.
  
\subsection{\bf Geometric flow}\hfill \break
\indent Let $M$ be a smooth, closed, uniformly convex hypersurface in $\mathbb{R}^n$ enclosing the origin in its interior, and $M$ be parametrised by the inverse spherical image $X=\nu_M^{-1}: \mathbb{S}^{n-1}\rightarrow M$. The support function $h$ of $M$ can be computed by
\begin{eqnarray}\label{3.1}
	h(x)=\langle x, X(x)\rangle,
\end{eqnarray}
where $x\in \mathbb{S}^{n-1}$ is the unit outer normal of $M$ at $X(x)$. It is easy to check that
\begin{eqnarray}\label{3.2}
	X(x)=h(x)x+\nabla h(x).
\end{eqnarray}
Here $\nabla$ is the covariant derivative with respect to the standard metric $e_{ij}$ of $\mathbb{S}^{n-1}$. The second fundamental form of $M$ is given by (see e.g. \cite{36a, 46a})
\begin{eqnarray*}
	H_{ij}=\nabla_{ij}h+he_{ij},
\end{eqnarray*}
where $\nabla_{ij}=\nabla^2$ denotes the second order covariant derivative with respect to $e_{ij}$. When we consider a smooth local orthonormal frame on $\mathbb{S}^{n-1}$, the principal radii of curvature of $M$ are the eigenvalues of the matrix
\begin{eqnarray*}
	b_{ij}=\nabla_{ij}h+h\delta_{ij}.
\end{eqnarray*}
Hence, the product of the principal radii of curvature is given by
\begin{eqnarray*}
	\mathcal{S}:=\det(\nabla_{ij}h+h\delta_{ij}).
\end{eqnarray*}
Let $u, x\in\mathbb{S}^{n-1}$ satisfy
\begin{eqnarray*}
	\rho(u)u=X(x)=h(x)x+\nabla h(x).
\end{eqnarray*}
Then the following formulae are well-known:
\begin{eqnarray}\label{3.3}
	\rho^2=h^2+|\nabla h|^2,~~h=\frac{\rho^2}{\sqrt{\rho^2+|\nabla \rho|^2}}.
\end{eqnarray}

Let $B_r$ be the smooth Euclidean ball, centered at the origin $o$, of radius $r>0$ in $\mathbb{R}^n$ and $\{M_t\}$ be a family of closed hypersufaces such that $M_t=X(\mathbb{S}^{n-1}, t)$. Taking the initial value $M_0=\partial B_r$ with $\gamma_n(B_r)=\frac{1}{2}$, and noting that $X:\mathbb{S}^{n-1}\times [0, T)\rightarrow \mathbb{R}^n$, we consider the following curvature flow
\begin{equation}\label{3.4}
	\left\{
	\begin{aligned}
		\frac{\partial X}{\partial t}(x, t)&=f(\nu)\mathcal{S}(x, t)e^{-\frac{\left|X\right|^2}{2}}\langle X, \nu\rangle\nu-( \sqrt{2\pi})^n \tau(t)X(x, t),\\
		X(x, 0)&=X_0(x)\in\partial B_r.
	\end{aligned}
	\right.
\end{equation}
Here, $f$ is a given positive smooth function on $\mathbb{S}^{n-1}$, $\nu$ is the unit outer normal vector of the hypersurface $M_t$ at the point $X(\cdot, t)$, $T$ is the maximal time for which the solution exists, $\mathcal{S}(\cdot, t)$ is the product of the principal radii of curvature with respect to $M_t$, and
\begin{eqnarray*}
	\tau(t)=\frac{\int_{\mathbb{S}^{n-1}}e^{-\frac{|X|^2}{2}}|X|^ndx}{( \sqrt{2\pi})^n\int_{\mathbb{S}^{n-1}}\frac{\langle X, \nu\rangle}{f}dx}.
\end{eqnarray*}
From the evolution equation (\ref{3.4}) of $X(x, t)$ and formula (\ref{3.1}), we see that the evolution equation of the corresponding support function $h(x, t)$ is 
\begin{equation}\label{3.5}
	\left\{
	\begin{aligned}
		\frac{\partial h}{\partial t}(x, t)&=f(x)\mathcal{S}(x, t)e^{-\frac{\rho(u, t)^2}{2}}h(x, t)-(\sqrt{2\pi})^n\tau(t)h(x, t),\\
		h(x, 0)&=r.
	\end{aligned}
	\right.
\end{equation}
Note that for $u\in \mathbb{S}^{n-1}$,  $\rho(u, t)$ is the radial function of $M_t$ with
\begin{eqnarray}\label{3.5+}
	\rho(u, t)=\sqrt{h(x, t)^2+|\nabla h(x, t)|^2}.
\end{eqnarray}
Moreover, it is easy to derive that
\begin{eqnarray}\label{3.6}
	\frac{1}{\rho(u, t)}\frac{\partial \rho(u, t)}{\partial t}=\frac{1}{h(x, t)}\frac{\partial h(x, t)}{\partial t}.
\end{eqnarray}
Then it follows from (\ref{3.5}) that the evolution equation of the radial function $\rho(u, t)$ is:
\begin{equation}\label{3.7}
	\left\{
	\begin{aligned}
		\frac{\partial \rho}{\partial t}(u, t)&=f(x)\mathcal{S}(x, t)e^{-\frac{\rho(u, t)^2}{2}}\rho(u, t)-(\sqrt{2\pi})^n\tau(t)\rho(u, t),\\
		\rho(u, 0)&=r,
	\end{aligned}
	\right.
\end{equation}
where $x=x(u, t)$ is the unit outer normal vector of $M_t$ at the point $\rho(u, t)u$.

Let $M_t$ be the convex hypersurface with the support function $h:\mathbb{S}^{n-1}\times [0, T)\rightarrow\mathbb{R}$. We define
\begin{eqnarray}\label{3.9}
	\mathcal{F}(t)=\frac{1}{|\mu|}\int_{\mathbb{S}^{n-1}}h(x, t)d\mu(x),
\end{eqnarray}
where $d\mu=\frac{1}{f}dx$ and $|\mu|=\mu\left(\mathbb{S}^{n-1}\right)$ for a Borel measure $\mu$ on $\mathbb{S}^{n-1}$. In the following, we see that the normalization factor $\frac{\int_{\mathbb{S}^{n-1}}e^{-\frac{\rho^2}{2}}\rho^ndu}{( \sqrt{2\pi})^n\int_{\mathbb{S}^{n-1}}h/fdx}$ garantees that the function (\ref{3.9}) remains unchanged under the flow (\ref{3.5}).

{\it \noindent{\bf Lemma 3.1.}~Along the flow (\ref{3.5}), we have
	\begin{eqnarray*}
		\mathcal{F}(0)=\mathcal{F}(t).
\end{eqnarray*}}
{\bf Proof.}~Note that
\begin{eqnarray*}
	\tau(t)=\frac{\int_{\mathbb{S}^{n-1}}e^{-\frac{\rho^2}{2}}\rho^ndu}{( \sqrt{2\pi})^n\int_{\mathbb{S}^{n-1}}h/fdx}
\end{eqnarray*}
and
\begin{eqnarray}\label{3.10}
	\frac{dx}{du}=\frac{\rho^n}{h\mathcal{S}}.
\end{eqnarray}
Then it follows from the flow (\ref{3.5}) and definition (\ref{3.9}) that
\begin{eqnarray*}
	\frac{d}{dt}\mathcal{F}(t)&=&\frac{1}{|\mu|}\int_{\mathbb{S}^{n-1}}\partial_th(x, t)\frac{1}{f}dx\\
	&=&\frac{1}{|\mu|}\int_{\mathbb{S}^{n-1}}\left( f\mathcal{S}he^{-\frac{\rho^2}{2}}-(\sqrt{2\pi})^n\tau(t)h\right) \frac{1}{f}dx\\
	&=&\frac{1}{|\mu|}\left(\int_{\mathbb{S}^{n-1}}\mathcal{S}
	he^{-\frac{\rho^2}{2}}dx-\frac{\int_{\mathbb{S}^{n-1}}e^{-\frac{\rho^2}{2}}\rho^ndu}{\int_{\mathbb{S}^{n-1}}h/fdx}\int_{\mathbb{S}^{n-1}}\frac{h}{f}dx\right)\\
	&=&0.
\end{eqnarray*}
This completes the proof of Lemma 3.1.\hfill $\square$

Let $K_t$ be the convex body enclosed by $M_t$. Then by the definition (\ref{1.0}) of Guassian volume, and using polar coordinates, we have
\begin{eqnarray}\label{3.11}
\gamma_n(t):=\gamma_n(K_t)=\frac{1}{(\sqrt{2\pi})^n}\int_{\mathbb{S}^{n-1}}\int_0^{\rho(u, t)}e^{-\frac{r^2}{2}}r^{n-1}drdu.
\end{eqnarray}

The following lemma shows that the function $\gamma_n(t)$ is non-decreasing along the flow (\ref{3.7}).

{\it \noindent{\bf Lemma 3.2.}~The function $\gamma_{n}(t)$ is non-decreasing along the flow (\ref{3.7}). That is,
	\begin{eqnarray*}
		\frac{d}{dt}\gamma_n(t)\geq 0,
	\end{eqnarray*}
and the equality holds if and only if $K_t$ satisfies the equation 
	\begin{eqnarray}\label{3.12}
	fe^{-\frac{\rho^2}{2}}\mathcal{S}=(\sqrt{2\pi})^n\tau(t).
\end{eqnarray}
Here, $\tau(t)=\int_{\mathbb{S}^{n-1}}e^{-\frac{\rho^2}{2}}\rho^ndu/(\sqrt{2\pi})^n\int_{\mathbb{S}^{n-1}}\frac{h}{f}dx$.}

{\bf Proof.}~From definition (\ref{3.11}) and flow (\ref{3.7}), we calculate
\begin{eqnarray*}
	&&(\sqrt{2\pi})^n\frac{d}{dt}\gamma_n(t)=\int_{\mathbb{S}^{n-1}}e^{-\frac{\rho^2}{2}}\rho^{n-1}\partial_t \rho(u, t)du \\
	&=&\int_{\mathbb{S}^{n-1}}e^{-\frac{\rho^2}{2}}\rho^{n-1}\left(f\mathcal{S}\rho e^{-\frac{\rho^2}{2}}-(\sqrt{2\pi})^n\tau(t)\rho\right)du \\
	&=&\int_{\mathbb{S}^{n-1}}f\mathcal{S}e^{-\rho^2}\rho^{n}du-\frac{\int_{\mathbb{S}^{n-1}}e^{-\frac{\rho^2}{2}}\rho^ndu}{\int_{\mathbb{S}^{n-1}}\frac{h}{f}dx}\int_{\mathbb{S}^{n-1}}e^{-\frac{\rho^2}{2}}\rho^ndu.
\end{eqnarray*}
Using the fact (\ref{3.10}), we have
\begin{eqnarray*}
	&&(\sqrt{2\pi})^n\left(\int_{\mathbb{S}^{n-1}}\frac{h}{f}dx\right)\gamma'_n(t) \\
	&=&\int_{\mathbb{S}^{n-1}}f\mathcal{S}e^{-\rho^2}\rho^{n}du\int_{\mathbb{S}^{n-1}}\frac{h}{f}dx-\int_{\mathbb{S}^{n-1}}e^{-\frac{\rho^2}{2}}\rho^ndu\int_{\mathbb{S}^{n-1}}e^{-\frac{\rho^2}{2}}\rho^ndu\\
	&=&\int_{\mathbb{S}^{n-1}}\left( \sqrt{f\mathcal{S}e^{-\rho^2}\rho^{n}}\right)^2 du\int_{\mathbb{S}^{n-1}}\left(\sqrt{\frac{\rho^n}{f\mathcal{S}}}\right)^2du-\left(\int_{\mathbb{S}^{n-1}}e^{-\frac{\rho^2}{2}}\rho^ndu\right)^2\\
	&\geq& \left(\int_{\mathbb{S}^{n-1}}e^{-\frac{\rho^2}{2}}\rho^ndu\right)^2-\left(\int_{\mathbb{S}^{n-1}}e^{-\frac{\rho^2}{2}}\rho^ndu\right)^2\\
	&=&0,
\end{eqnarray*}
where the above inequality is due to the H\"{o}lder inequality. From the equation condition of H\"{o}lder inequality, we see that the equality of the inequality holds if and only if 
	\begin{eqnarray*}
	\frac{1}{(\sqrt{2\pi})^n}fe^{-\frac{\rho^2}{2}}\mathcal{S}=\tau(t)~~\mbox{with}~~\tau(t)=\frac{\int_{\mathbb{S}^{n-1}}e^{-\frac{\rho^2}{2}}\rho^ndu}{( \sqrt{2\pi})^n\int_{\mathbb{S}^{n-1}}h/fdx}.
\end{eqnarray*}
The proof is completed. \hfill $\square$

\subsection{\bf Long-time existence of the flow}\hfill \break
\indent In this subsection, we establish the priori estimates and show that the flow (\ref{3.4}) exists for all time. We begin with the estimates of the uniformly positive lower and upper bounds of $h, \rho, \tau$. The following function will be required.
\begin{eqnarray}\label{3.13}
	\left\| g\right\|_\mu=\frac{1}{|\mu|}\int_{\mathbb{S}^{n-1}}g(x)d\mu(x)
\end{eqnarray}
for a continuous $g: \mathbb{S}^{n-1}\rightarrow [0, \infty)$.

{\it \noindent{\bf Lemma 3.3.}~Let $X(\cdot, t)$ be a uniformly convex solution to the flow (\ref{3.4}). For a positive and smooth function $f$ on $\mathbb{S}^{n-1}$, suppose that $\mu$ is a finite Borel measure not concentrated in any closed hemisphere with $d\mu=\frac{1}{f}dx$. Then there exists a positive constant $C$ independent of $t$ such that for every $t\in [0, T)$,
	\begin{eqnarray}\label{3.14}
		\frac{1}{C}\leq h(\cdot, t), ~~\rho(\cdot, t)\leq C~\mbox{on}~~\mathbb{S}^{n-1},
	\end{eqnarray}
			\begin{eqnarray*}
		\frac{1}{C}\leq \tau(t)\leq C.
	\end{eqnarray*}
	}
{\bf Proof.}~From (\ref{2.4}), we only need to establish the upper bound of $h(\cdot, t)$ and the lower bound of $\rho(\cdot, t)$ for inequalities (\ref{3.14}).

Let $u\in \mathbb{S}^{n-1}$ and $R_t>0$ such that $R_t u$ is maximal for $R_tu\in M_t=X(\mathbb{S}^{n-1}, t)$. This means 
\begin{equation}\label{3.16}
R_th_{\tilde{u}}\leq h(\cdot, t).
\end{equation}
From Lemma 3.1, we see
\begin{equation}\label{3.17}
	\mathcal{F}(t)=	\mathcal{F}(0)=r>0.
\end{equation}
Thus from (\ref{3.17}) we have
\begin{equation}\label{3.18}
\frac{1}{|\mu|}\int_{\mathbb{S}^{n-1}}h(x, t)d\mu(x)=\mathcal{F}(t)=r.
\end{equation}
Combining (\ref{3.16}) and (\ref{3.18}), it follows that
\begin{equation*}
	r\geq R_t\left\|h_{\tilde{u}}\right\|_\mu\geq R_t\min_{u\in \mathbb{S}^{n-1}}\left\|h_{\tilde{u}}\right\|_\mu.
\end{equation*}
Namely,
\begin{equation}\label{3.19}
h(\cdot, t)\leq R_t\leq \frac{r}{\min_{u\in \mathbb{S}^{n-1}}\left\|h_{\tilde{u}}\right\|_\mu}.
\end{equation}
Further, we show that there exists a constant $c>0$ such that 
\begin{equation}\label{3.20}
\min_{u\in \mathbb{S}^{n-1}}\left\|h_{\tilde{u}}\right\|_\mu\geq c>0.
\end{equation}
Since the support of $\mu$ is not contained in a closed hemisphere of $\mathbb{S}^{n-1}$, we have that for every $u\in \mathbb{S}^{n-1}$, 
\begin{equation}\label{3.21}
\mu\left(\left\{h_{\tilde{u}}>0\right\}\right)=\mu\left(\mathbb{S}^{n-1}\backslash H^-_u\right)>0,
\end{equation}
where
\begin{equation*}
	H^-_u:=\left\{y\in \mathbb{R}^n: \langle u, y \rangle\leq 0\right\}
\end{equation*}
denotes the closed halfspaces with outer normal vector $u$. Thus, by the fact (\ref{3.21}) and definition (\ref{3.13}) we have
\begin{equation*}
\left\|h_{\tilde{u}}\right\|_\mu=\frac{1}{|\mu|}\int_{\mathbb{S}^{n-1}}h_{\tilde{u}}(x)d\mu(x)\geq 	\frac{1}{|\mu|}\int_{\{h_{\tilde{u}}>0\}}h_{\tilde{u}}(x)d\mu(x)>0.
\end{equation*}
Let $u_i\in \mathbb{S}^{n-1}$ with $u_i\rightarrow u \in\mathbb{S}^{n-1}$ as $i\rightarrow +\infty$. This implies $h_{\tilde{u}_i}\rightarrow h_{\tilde{u}}$ uniformly on $\mathbb{S}^{n-1}$. Note that $h_{\tilde{u}_i}\leq 1$. Then from the dominated convergence theorem we have
\begin{eqnarray*}
\lim_{i\rightarrow \infty}\left\|h_{\tilde{u}_i} \right\|_\mu = \lim_{i\rightarrow \infty}\frac{1}{|\mu|}\int_{\mathbb{S}^{n-1}}h_{\tilde{u}_i}(x)d\mu(x) =\frac{1}{|\mu|}\int_{\mathbb{S}^{n-1}}h_{\tilde{u}}(x)d\mu(x)=\left\|h_{\tilde{u}} \right\|_\mu. 
\end{eqnarray*}
That is, the function $u\mapsto \left\|h_{\tilde{u}} \right\|_\mu$ is continuous. Noting that the unit sphere $\mathbb{S}^{n-1}$ is compact, it follows that there exists a constant $c>0$ such that 
\begin{eqnarray*}
\min_{u\in \mathbb{S}^{n-1}}\left\|h_{\tilde{u}} \right\|_\mu\geq c>0. 
\end{eqnarray*}
Recalling (\ref{3.19}), we see that $h(\cdot, t)$ has the uniformly positive upper bound independent of $t$. 

For the lower bound of $\rho(\cdot, t)$, we prove it by contradiction. Let $\rho(\cdot, t)\rightarrow 0$. Note that $\rho(\cdot, t)$ has the uniformly upper bound from (\ref{2.4}). Thus it follows from the dominated convergence theorem that
\begin{eqnarray*}
\gamma_n(t)=\frac{1}{(\sqrt{2\pi})^n}\int_{\mathbb{S}^{n-1}}\int_0^{\rho(u, t)}e^{-\frac{r^2}{2}}r^{n-1}drdu\rightarrow 0.
\end{eqnarray*}
This is a contradiction since
\begin{eqnarray*}
\gamma_n(t)\geq \gamma_n(0)=\frac{1}{2}
\end{eqnarray*}
from Lemma 3.2. Hence $\rho(\cdot, t)$ is uniformly bounded from below.

Since 
\begin{eqnarray*}
\tau(t)=\frac{\int_{\mathbb{S}^{n-1}}e^{-\frac{\rho^2}{2}}\rho^ndu}{( \sqrt{2\pi})^n\int_{\mathbb{S}^{n-1}}h/fdx},
\end{eqnarray*}
the estimates about $\tau(t)$ follows directly from the lower and upper bounds of $h(\cdot, t)$ and $\rho(\cdot, t)$.
\hfill $\square$

From the identities in (\ref{3.3}) and the inequalities in (\ref{3.14}), we have the following gradient estimates of $h(\cdot, t)$ and $\rho(\cdot, t)$.

{\it \noindent{\bf Corollary 3.1.}~There exists a positive constant $C$ independent of $t$ such that for every $t\in [0, T)$,
	\begin{eqnarray*}
		|\nabla h(\cdot, t)|\leq C~\mbox{on}~~\mathbb{S}^{n-1},
	\end{eqnarray*}
	\begin{eqnarray*}
		|\nabla \rho(\cdot, t)|\leq C~\mbox{on}~~\mathbb{S}^{n-1}.
\end{eqnarray*}}
In order to obtain the long-time existence of the flow (\ref{3.4}), we further need to establish the lower and upper bounds of the eigenvalues of the matrix $b_{ij}$, which are given in the next Lemma 3.5. We first derive Lemma 3.4, i.e., the uniformly upper bound for the product of the principal radii of curvature will be eatablished. For notational simplicity, we write
\begin{eqnarray*}
	\mathcal{Z}=fe^{-\frac{\rho^2}{2}}h.
\end{eqnarray*} 

In the rest of this section, if in a term the same index appears twice (both as lower and upper index), then the term is assumed to be summed from $1$ to $n-1$. For example,
	\begin{eqnarray*}
a_ib^i=\sum_{i=1}^{n-1}a_ib^i, ~~~~~~~~ a^{ij}b_{ij}=\sum_{i,j=1}^{n-1}a^{ij}b_{ij}.
\end{eqnarray*}

{\it \noindent{\bf Lemma 3.4.}~There is a positive constant $C$ independent of $t$ such that for every $t\in [0, T)$,
	\begin{eqnarray*}
		\mathcal{S}(\cdot, t)\leq C~\mbox{on}~~\mathbb{S}^{n-1}.
\end{eqnarray*}}
{\bf Proof.}~Let $\varepsilon>0$ be small enough. We consider the following auxiliary function
\begin{eqnarray}\label{3.21+}
	P=\frac{\frac{\mathcal{Z}\mathcal{S}}{h}}{1-\varepsilon \frac{\rho^2}{2}}.
\end{eqnarray}
We write $b^{ij}$ for the inverse matrix of $b_{ij}$. Then the eigenvalues of $b_{ij}$ and $b^{ij}$ are respectively the principal radii and principal curvatures of $M_t$. For any fixed $t\in [0, T)$, assume that $\max_{\mathbb{S}^{n-1}}P(x, t)$ is attained at some point on $\mathbb{S}^{n-1}$. Denote by $\mathcal{S}^{ij}$ the cofactor matrix of $b_{ij}$. If we take an orthonormal frame such that $b_{ij}$ and $b^{ij}$ are diagonal at this point, then we first compute the evolution equation of $\frac{\mathcal{Z}\mathcal{S}}{h}$.
\begin{eqnarray*}
	\nabla_i\left(\frac{\mathcal{Z}\mathcal{S}}{h}\right)=\frac{\nabla_i(\mathcal{Z}\mathcal{S})
	}{h}-\frac{\mathcal{Z}\mathcal{S}\nabla_ih}{h^2},
\end{eqnarray*}

\begin{eqnarray}\label{3.22}
	\nabla_{ij}\left(\frac{\mathcal{Z}\mathcal{S}}{h}\right)
	&=&\frac{\nabla_{ij}(\mathcal{Z}\mathcal{S})}{h}-\frac{\mathcal{Z}\mathcal{S}\nabla_{ij}h}{h^2}-\frac{2h\nabla_ih\nabla_j(\mathcal{Z}
		\mathcal{S})-2\mathcal{Z}\mathcal{S}\nabla_ih\nabla_jh}{h^3}\nonumber\\
	&=&\frac{\nabla_{ij}(\mathcal{Z}\mathcal{S})}{h}-\frac{\mathcal{Z}\mathcal{S}\nabla_{ij}h}{h^2}-2\nabla_i(\log h)\nabla_j\left(\frac{\mathcal{Z}\mathcal{S}}{h}\right),
\end{eqnarray}
\begin{eqnarray}\label{3.23}
	\partial_t(\mathcal{Z}\mathcal{S})&=&\mathcal{Z}\mathcal{S}^{ij}(\nabla_{ij}(\partial_th)+\partial_th\delta_{ij})+\mathcal{S}\partial_t\mathcal{Z}\nonumber \\
	&=&\mathcal{Z}\mathcal{S}^{ij}(\nabla_{ij}(\partial_th)+\partial_th\delta_{ij})
	+\mathcal{S}f\partial_t\left(e^{-\frac{\rho^2}{2}}h\right)\nonumber \\
	&=&\mathcal{Z}\mathcal{S}^{ij}(\nabla_{ij}(\partial_th)+\partial_th\delta_{ij})-\mathcal{S}f\rho h e^{-\frac{\rho^2}{2}}\partial_t\rho+\mathcal{S}f e^{-\frac{\rho^2}{2}}\partial_th \nonumber \\
	&=&\mathcal{Z}\mathcal{S}^{ij}\left[\nabla_{ij}(\mathcal{Z}\mathcal{S}-(\sqrt{2\pi})^n\tau h)+(\mathcal{Z}\mathcal{S}-(\sqrt{2\pi})^n\tau h)\delta_{ij}\right]\nonumber \\
	&&-\mathcal{S}f\rho h e^{-\frac{\rho^2}{2}}\left(\frac{\rho}{h}\mathcal{Z}\mathcal{S}-(\sqrt{2\pi})^n\tau\rho\right)+\mathcal{S}f e^{-\frac{\rho^2}{2}}\left(\mathcal{Z}\mathcal{S}-(\sqrt{2\pi})^n\tau h\right) \nonumber \\
	&=&\mathcal{Z}\mathcal{S}^{ij}\nabla_{ij}(\mathcal{Z}\mathcal{S})+\mathcal{S}^{ij}\mathcal{Z}^2\mathcal{S}\delta_{ij}-(\sqrt{2\pi})^n\tau\mathcal{Z}\mathcal{S}^{ij}
	(\nabla_{ij}h+h\delta_{ij}) \nonumber \\
	&&-\rho\mathcal{Z}\mathcal{S}\left(\frac{\rho}{h}\mathcal{Z}\mathcal{S}-(\sqrt{2\pi})^n\tau\rho\right)+\frac{\mathcal{Z}\mathcal{S}}{h}\left(\mathcal{Z}\mathcal{S}-(\sqrt{2\pi})^n\tau h\right) \nonumber \\
	&=&\frac{(\mathcal{Z}\mathcal{S})^2}{h}-\frac{\rho^2(\mathcal{Z}\mathcal{S})^2}{h}+\mathcal{S}^{ij}\mathcal{Z}^2\mathcal{S}\delta_{ij}+(\sqrt{2\pi})^n\tau\rho^2(\mathcal{Z}\mathcal{S})-(\sqrt{2\pi})^n\tau(\mathcal{Z}\mathcal{S})\nonumber \\
	&&-(n-1)(\sqrt{2\pi})^n\tau(\mathcal{Z}\mathcal{S})+\mathcal{Z}\mathcal{S}^{ij}\nabla_{ij}(\mathcal{Z}\mathcal{S})\nonumber \\
	&=&\left(\frac{1}{h}-\frac{\rho^2}{h}\right)(\mathcal{Z}\mathcal{S})^2
	+\left(\mathcal{S}^{ij}\delta_{ij}\right)\left(\mathcal{Z}^2\mathcal{S}\right)+(\sqrt{2\pi})^n\tau(\rho^2-n)(\mathcal{Z}\mathcal{S})\nonumber\\
	&&+\mathcal{Z}\mathcal{S}^{ij}\nabla_{ij}(\mathcal{Z}\mathcal{S}).
\end{eqnarray}
By (\ref{3.23}), we have
\begin{eqnarray}\label{3.24}
	\partial_t\left(\frac{\mathcal{Z}\mathcal{S}}{h}\right)
	&=&\frac{\partial_t(\mathcal{Z}\mathcal{S})}{h}-\frac{\mathcal{Z}\mathcal{S}}{h^2}\partial_th\nonumber\\
	&=&\left(\frac{1}{h^2}-\frac{\rho^2}{h^2}\right)(\mathcal{Z}\mathcal{S})^2
	+\frac{1}{h}\left(\mathcal{S}^{ij}\delta_{ij}\right)\left(\mathcal{Z}^2\mathcal{S}\right)+(\sqrt{2\pi})^n\tau\left(\frac{\rho^2-n}{h}\right)(\mathcal{Z}\mathcal{S})\nonumber\\
	&&
	+\frac{1}{h}\mathcal{Z}\mathcal{S}^{ij}\nabla_{ij}(\mathcal{Z}\mathcal{S})-\frac{\mathcal{Z}\mathcal{S}}{h^2}\left(\mathcal{Z}\mathcal{S}-(\sqrt{2\pi})^n\tau h\right).
\end{eqnarray}
Then according to (\ref{3.22}) and (\ref{3.24}), the evolution equation of $\frac{\mathcal{Z}\mathcal{S}}{h}$ is given by
\begin{eqnarray*}
	&&\partial_t\left(\frac{\mathcal{Z}\mathcal{S}}{h}\right)-\mathcal{Z}\mathcal{S}^{ij}\nabla_{ij}\left(\frac{\mathcal{Z}\mathcal{S}}{h}\right)\nonumber\\
	&=&\left(\frac{1}{h^2}-\frac{\rho^2}{h^2}\right)(\mathcal{Z}\mathcal{S})^2
	+\frac{1}{h}\left(\mathcal{S}^{ij}\delta_{ij}\right)\left(\mathcal{Z}^2\mathcal{S}\right)-\frac{1}{h^2}(\mathcal{Z}\mathcal{S})^2+\frac{\mathcal{S}^{ij}\nabla_{ij}h}{h^2}\left(\mathcal{Z}^2\mathcal{S}\right)\nonumber\\
	&&	-\frac{(n-1)(\sqrt{2\pi})^n\tau}{h}\left(\mathcal{Z}\mathcal{S}\right)+\frac{(\sqrt{2\pi})^n\tau\rho^2}{h}\left(\mathcal{Z}\mathcal{S}\right)+2\mathcal{S}^{ij}\nabla_i(\log h)\nabla_j\left(\frac{\mathcal{Z}\mathcal{S}}{h}\right)\mathcal{Z}.\nonumber\\
&=&	2\mathcal{S}^{ij}\nabla_i(\log h)\nabla_j\left(\frac{\mathcal{Z}\mathcal{S}}{h}\right)\mathcal{Z}-\frac{\rho^2}{h^2}(\mathcal{Z}\mathcal{S})^2+\frac{1}{h}\left(\mathcal{S}^{ij}\delta_{ij}\right)\left(\mathcal{Z}^2\mathcal{S}\right)\nonumber\\
&&+\frac{(\sqrt{2\pi})^n\tau\left(\rho^2-n+1\right)}{h}\left(\mathcal{Z}\mathcal{S}\right)+\frac{\mathcal{S}^{ij}\left(b_{ij}-h\delta_{ij}\right)}{h^2}\left(\mathcal{Z}^2\mathcal{S}\right)\nonumber\\
\end{eqnarray*}
\begin{eqnarray}\label{3.25}	
&=&-\frac{\left(\rho^2-n+1\right)}{h^2}(\mathcal{Z}\mathcal{S})^2+\frac{(\sqrt{2\pi})^n\tau\left(\rho^2-n+1\right)}{h}\left(\mathcal{Z}\mathcal{S}\right) \nonumber\\	
&&+2\mathcal{S}^{ij}\nabla_i(\log h)\nabla_j\left(\frac{\mathcal{Z}\mathcal{S}}{h}\right)\mathcal{Z}.
\end{eqnarray}
Now, we compute the evolution equation of $P$.
\begin{eqnarray}\label{3.26}
	0=\nabla_i\left(\frac{\frac{\mathcal{Z}\mathcal{S}}{h}}{1-\varepsilon \frac{\rho^2}{2}}\right)	&=&\frac{\nabla_i\left(\frac{\mathcal{Z}\mathcal{S}}{h}\right)}{1-\varepsilon \frac{\rho^2}{2}}+\frac{\varepsilon\frac{\mathcal{Z}\mathcal{S}}{h}\nabla_i\left
		(\frac{\rho^2}{2}\right)}{\left(1-\varepsilon \frac{\rho^2}{2}\right)^2}.
\end{eqnarray}
Similarly, 
\begin{eqnarray}\label{3.27}
	0=\nabla_j\left(\frac{\frac{\mathcal{Z}\mathcal{S}}{h}}{1-\varepsilon \frac{\rho^2}{2}}\right)=\frac{\nabla_j\left(\frac{\mathcal{Z}\mathcal{S}}{h}\right)}{1-\varepsilon \frac{\rho^2}{2}}+\frac{\varepsilon\frac{\mathcal{Z}\mathcal{S}}{h}\nabla_j\left
		(\frac{\rho^2}{2}\right)}{\left(1-\varepsilon \frac{\rho^2}{2}\right)^2}.
\end{eqnarray} 
Together (\ref{3.26}) with (\ref{3.27}), it follows that
 \begin{eqnarray}\label{3.28}
  \ \ \ \ \ \  \ \ \	\frac{\varepsilon\nabla_i\left(\frac{\mathcal{Z}\mathcal{S}}{h}\right)\nabla_j\left
 		(\frac{\rho^2}{2}\right)}{\left(1-\varepsilon \frac{\rho^2}{2}\right)^2}+\frac{\varepsilon\nabla_j\left(\frac{\mathcal{Z}\mathcal{S}}{h}\right)\nabla_i\left
 		(\frac{\rho^2}{2}\right)}{\left(1-\varepsilon \frac{\rho^2}{2}\right)^2}=-\frac{2\varepsilon^2\left(\frac{\mathcal{Z}\mathcal{S}}{h}\right)\nabla_i\left(\frac{\rho^2}{2}\right)\nabla_j\left
 		(\frac{\rho^2}{2}\right)}{\left(1-\varepsilon \frac{\rho^2}{2}\right)^3}.
 \end{eqnarray}
Then by (\ref{3.28}) we have
 \begin{eqnarray}\label{3.29}
 	\nabla_{ij}\left(\frac{\frac{\mathcal{Z}\mathcal{S}}{h}}{1-\varepsilon \frac{\rho^2}{2}}\right)
 	&=&\frac{\nabla_{ij}\left(\frac{\mathcal{Z}\mathcal{S}}{h}\right)}{1-\varepsilon \frac{\rho^2}{2}}+\frac{\varepsilon\nabla_i\left(\frac{\mathcal{Z}\mathcal{S}}{h}\right)\nabla_j\left
 		(\frac{\rho^2}{2}\right)}{\left(1-\varepsilon \frac{\rho^2}{2}\right)^2}+\frac{\varepsilon\nabla_i\left
 		(\frac{\rho^2}{2}\right)\nabla_j\left(\frac{\mathcal{Z}\mathcal{S}}{h}\right)}{\left(1-\varepsilon \frac{\rho^2}{2}\right)^2}\nonumber\\
 	&&+\frac{\varepsilon\frac{\mathcal{Z}\mathcal{S}}{h} \nabla_{ij}\left
 		(\frac{\rho^2}{2}\right)}{\left(1-\varepsilon\frac{\rho^2}{2}\right)^2}+\frac{2\varepsilon^2\frac{\mathcal{Z}\mathcal{S}}{h}\nabla_i\left
 		(\frac{\rho^2}{2}\right)\nabla_j\left
 		(\frac{\rho^2}{2}\right)}{\left(1-\varepsilon \frac{\rho^2}{2}\right)^3}\nonumber\\
 	&=&\frac{\nabla_{ij}\left(\frac{\mathcal{Z}\mathcal{S}}{h}\right)}{1-\varepsilon \frac{\rho^2}{2}}+\frac{\varepsilon\frac{\mathcal{Z}\mathcal{S}}{h}\nabla_{ij}\left
 		(\frac{\rho^2}{2}\right)}{\left(1-\varepsilon \frac{\rho^2}{2}\right)^2}.
 \end{eqnarray}
 Recalling (\ref{3.5+}), we obtain
 \begin{eqnarray*}
 	\nabla_i\left(\frac{\rho^2}{2}\right)&=&h\nabla_ih+\sum_k\nabla_kh\nabla_{ki}h.
 \end{eqnarray*}
By this, it follows that
 \begin{eqnarray}\label{3.30}
 	\nabla_{ij}\left(\frac{\rho^2}{2}\right)&=&\nabla_ih\nabla_jh+h\nabla_{ij}h+\sum_k\nabla_{ki}h\nabla_{kj}h+\sum_k\nabla_{k}h\nabla_{kij}h \nonumber \\
 	&=&\nabla_ih\nabla_jh+h(b_{ij}-h\delta_{ij})+\sum_k\nabla_kh\nabla_j(b_{ki}-h\delta_{ki})\nonumber \\
 	&&+\sum_k(b_{ki}-h\delta_{ki})(b_{kj}-h\delta_{kj})\nonumber \\	
 	&=&\nabla_ih\nabla_jh+hb_{ij}-h^2\delta_{ij}+\sum_k\nabla_kh\nabla_jb_{ki}-\sum_k\nabla_kh\nabla_jh\delta_{ki}\nonumber \\
 	&&+\sum_kb_{ki}b_{kj}-h\sum_kb_{ki}\delta_{kj}-h\sum_kb_{kj}\delta_{ki}+h^2\sum_k\delta_{ki}\delta_{kj}.
 \end{eqnarray}
 Hence,
 \begin{eqnarray*}
 	\mathcal{S}^{ij}\nabla_{ij}\left(\frac{\rho^2}{2}\right)&=&\mathcal{S}^{ij}\nabla_ih\nabla_jh+(n-1)h\mathcal{S}-h^2\mathcal{S}^{ij}
 	\delta_{ij}+\sum_k\nabla_kh\nabla_k\mathcal{S}-(n-1)h\mathcal{S}\\
 	&&-\mathcal{S}^{ij}\sum_k\nabla_kh\nabla_jh\delta_{ki}+\mathcal{S}^{ij}\sum_kb_{ki}b_{kj}-(n-1)h\mathcal{S}+h^2\mathcal{S}^{ij}\delta_{ij}\\
 	&=&-(n-1)h\mathcal{S}+\sum_k\nabla_kh\nabla_k\mathcal{S}+\mathcal{S}^{ij}\sum_kb_{ki}b_{kj}.
 \end{eqnarray*}
 Moreover,
 \begin{eqnarray}\label{3.31}
 	\partial_t\left(\frac{\rho^2}{2}\right)&=&h\partial_th+\sum_k\nabla_kh\nabla_{k}(\partial_th)\nonumber \\
 	&=&h\left(\mathcal{Z}\mathcal{S}-(\sqrt{2\pi})^n\tau h\right)+\sum_k\nabla_kh\left(\mathcal{S}\nabla_k\mathcal{Z}+\mathcal{Z}\nabla_k\mathcal{S}-(\sqrt{2\pi})^n\tau\nabla_kh\right)\nonumber \\
 	&=&h\mathcal{Z}\mathcal{S}-(\sqrt{2\pi})^n\tau h^2+\sum_k(\nabla_kh\nabla_k\mathcal{Z})\mathcal{S}\nonumber \\
 &&+\sum_k(\nabla_kh\nabla_k\mathcal{S})\mathcal{Z}-(\sqrt{2\pi})^n\tau|\nabla h|^2\nonumber \\
 	&=&h\mathcal{Z}\mathcal{S}-(\sqrt{2\pi})^n\tau \rho^2+\mathcal{S}\sum_k(\nabla_kh\nabla_k\mathcal{Z})+\mathcal{Z}\sum_k(\nabla_kh\nabla_k\mathcal{S}).
 \end{eqnarray}
 From the above results, we have the evolution equation of $\frac{\rho^2}{2}$:
\begin{eqnarray}\label{3.32}
	&&\partial_t\left(\frac{\rho^2}{2}\right)-\mathcal{Z}\mathcal{S}^{ij}\nabla_{ij}\left(\frac{\rho^2}{2}\right)\nonumber \\
	&=&nh\mathcal{Z}\mathcal{S}
	+\mathcal{S}\sum_k(\nabla_kh\nabla_k\mathcal{Z})-\mathcal{Z}\mathcal{S}^{ij}\sum_kb_{ki}b_{kj}-(\sqrt{2\pi})^n\tau\rho^2.
\end{eqnarray} 
Making use of (\ref{3.32}), (\ref{3.25}) and (\ref{3.29}), we obtain the evolution equation of $P$: 
 \begin{eqnarray}\label{3.33}
 	&&\frac{\partial P}{\partial t}-\mathcal{Z}\mathcal{S}^{ij}\nabla_{ij}P\nonumber\\
 	&=&\frac{\partial_t\left(\frac{\mathcal{Z}\mathcal{S}}{h}\right)\left(1-\varepsilon \frac{\rho^2}{2}\right)+\varepsilon\frac{\mathcal{Z}\mathcal{S}}{h} \partial_t\left(\frac{\rho^2}{2}\right)}{\left(1-\varepsilon \frac{\rho^2}{2}\right)^2}-\mathcal{Z}\mathcal{S}^{ij}\nabla_{ij}P\nonumber\\	
  &=&\frac{\partial_t\left(\frac{\mathcal{Z}\mathcal{S}}{h}\right)}{1-\varepsilon \frac{\rho^2}{2}}+\frac{\varepsilon\frac{\mathcal{Z}\mathcal{S}}{h}\partial_t\left(\frac{\rho^2}{2}\right)}{\left(1-\varepsilon \frac{\rho^2}{2}\right)^2}-
 \frac{\mathcal{Z}\mathcal{S}^{ij}\nabla_{ij}\left(\frac{\mathcal{Z}\mathcal{S}}{h}\right)}{1-\varepsilon \frac{\rho^2}{2}}-\frac{\varepsilon\mathcal{Z}\mathcal{S}^{ij}\frac{\mathcal{Z}\mathcal{S}}{h}\nabla_{ij}\left(\frac{\rho^2}{2}\right)}{\left(1-\varepsilon \frac{\rho^2}{2}\right)^2}\nonumber\\
 &=&\frac{1}{1-\varepsilon \frac{\rho^2}{2}}\left[\partial_t\left(\frac{\mathcal{Z}\mathcal{S}}{h}\right)-\mathcal{Z}\mathcal{S}^{ij}
 \nabla_{ij}\left(\frac{\mathcal{Z}\mathcal{S}}{h}\right)\right]\nonumber\\
 &&+\frac{\varepsilon\frac{\mathcal{Z}\mathcal{S}}{h}}{\left(1-\varepsilon \frac{\rho^2}{2}\right)^2}\left[\partial_t\left(\frac{\rho^2}{2}\right)-\mathcal{Z}\mathcal{S}^{ij}\nabla_{ij}\left(\frac{\rho^2}{2}\right)\right]\nonumber\\
 &=&\frac{1}{1-\varepsilon \frac{\rho^2}{2}}\left[2\mathcal{S}^{ij}\nabla_i(\log h)\nabla_j\left(\frac{\mathcal{Z}\mathcal{S}}{h}\right)\mathcal{Z}-\frac{\left(\rho^2-n+1\right)}{h^2}(\mathcal{Z}\mathcal{S})^2\right.\nonumber\\
 &&\left.+\frac{(\sqrt{2\pi})^n\tau\left(\rho^2-n+1\right)}{h}\left(\mathcal{Z}\mathcal{S}\right) \right]\nonumber\\
 &&+ \frac{\varepsilon\frac{\mathcal{Z}\mathcal{S}}{h}}{\left(1-\varepsilon \frac{\rho^2}{2}\right)^2}\left[  nh\mathcal{Z}\mathcal{S}
 +\mathcal{S}\sum_k(\nabla_kh\nabla_k\mathcal{Z})-\mathcal{Z}\mathcal{S}^{ij}\sum_kb_{ki}b_{kj}-(\sqrt{2\pi})^n\tau\rho^2\right]. 
 \end{eqnarray}
 Recalling (\ref{3.27}), we have
\begin{eqnarray}\label{3.34}
	\nabla_j\left(\frac{\mathcal{Z}\mathcal{S}}{h}\right)=-\frac{\varepsilon\frac{\mathcal{Z}\mathcal{S}}{h}\nabla_j\left
		(\frac{\rho^2}{2}\right)}{1-\varepsilon \frac{\rho^2}{2}}.
\end{eqnarray} 
Note that 
\begin{eqnarray}\label{3.35}
	\nabla_j\left(\frac{\rho^2}{2}\right)=h\nabla_jh+\sum_k\nabla_kh\nabla_{kj}h=\sum_k\nabla_{k}hb_{kj}.
\end{eqnarray} 
Then from (\ref{3.34}) and (\ref{3.35}) we have
\begin{eqnarray}\label{3.36}
\  \ \ \ \ \ \  \ \ 	\mathcal{S}^{ij}\nabla_i(\log h)\nabla_j\left(\frac{\mathcal{Z}\mathcal{S}}{h}\right)=\frac{\mathcal{S}^{ij}}{h}\nabla_ih\nabla_j\left(\frac{\mathcal{Z}\mathcal{S}}{h}\right)
	=-\frac{\varepsilon(n-1)\mathcal{S}^2\mathcal{Z}|\nabla h|^2}{h^2\left(1-\varepsilon \frac{\rho^2}{2}\right)}\leq 0.
\end{eqnarray}
Using the arithmetic-geometric mean inequality, we give
\begin{eqnarray}\label{3.37}
	\mathcal{S}^{ij}\sum_kb_{ik}b_{jk}\geq (n-1)\mathcal{S}^{1+\frac{1}{n-1}}.
\end{eqnarray}
By Lemma 3.3, Corollary 3.1, (\ref{3.36}) and (\ref{3.37}), it follows from (\ref{3.33}) that
\begin{eqnarray*}
&&\frac{\partial P}{\partial t}-\mathcal{Z}\mathcal{S}^{ij}\nabla_{ij}P \\
&\leq& C_1 P+C_2 P^2-C_3 P^{2+\frac{1}{n-1}}\leq C_3P^2\left(C_4-P^{\frac{1}{n-1}}\right)
\end{eqnarray*}
for some positive constants $C_1, C_2, C_3, C_4$. This implies
 \begin{eqnarray*}
 \frac{\partial P}{\partial t}<0
 \end{eqnarray*}
for a large $P$. Thus, we obtain that for every $t\in [0, T)$,
\begin{eqnarray*}
	P(\cdot, t)\leq C~~\mbox{on}~ \mathbb{S}^{n-1}.
\end{eqnarray*}
By Lemma 3.3, and noting that (\ref{3.21+}), we have that for a positive constant $C_0$,
\begin{eqnarray*}
	\frac{1}{C_0}\mathcal{S}\leq P\leq C_0\mathcal{S}.
\end{eqnarray*}
This gives that $\mathcal{S}$ has a uniformly upper bound.\hfill $\square$
 
Next, we prove that the eigenvalues of the matrix $b_{ij}$ are bounded by positive constants from both above and below.
 
{\it \noindent{\bf Lemma 3.5.}~There is a positive constant $C$ independent of $t$ such that for $i=1,\cdots, n-1$, the principal radii $\lambda_i$ of curvatures of the hypersufaces $M_t$ satisfy
	\begin{eqnarray}\label{3.38}
		\frac{1}{C}\leq\lambda_i(\cdot, t)\leq C~\mbox{on}~~\mathbb{S}^{n-1}\times [0, T).
\end{eqnarray}} 
\ \ {\bf Proof.}~Let $M^*_t$ be the polar set of $M_t=X(\mathbb{S}^{n-1}, \cdot)$ and $h^*(\cdot, t)$ be the support function of $M_t^*$. If $\rho(\cdot, t)$ is the radial function of $M_t$, then from (\ref{2.3}) we have 
\begin{eqnarray}\label{3.39}
	\rho(\cdot, t)=\frac{1}{h^*(\cdot, t)}.
\end{eqnarray}
Let $p\in M_t$ and $p^*\in M^*_t$ satisfy $\langle p, p^*\rangle=1$. Then the following fact was given in \cite{19a, 44a}:
\begin{eqnarray}\label{3.40}
	h(x, t)^{n+1}h^*(u, t)^{n+1}\mathcal{S}(p)\mathcal{S}^*(p^*)=1,
\end{eqnarray}
where $x$ and $u$ are the unit outer normals of $M_t$ and $M^*_t$ at $p$ and $p^*$, respectively. Moreover, it follows from Lemma 3.3 and (\ref{2.3}) that for a positive constant $C$ independent of $t$,
\begin{eqnarray}\label{3.41}
\frac{1}{C}\leq h^*(\cdot, t)\leq C.
\end{eqnarray}
Hence, from Lemmas 3.3 and 3.4 as well as inequalities (\ref{3.41}), the above equation (\ref{3.40}) implies that there exists a positive constant $C$ independent of $t$ such that 
\begin{eqnarray}\label{3.42}
	\frac{1}{C}\leq \mathcal{S}^*(\cdot, t).
\end{eqnarray}
For $i=1, \cdots, n-1$, let $\lambda^*_i$ be the principal radii of curvatures of $M_t^*$. Then by the estimate (\ref{3.42}) and duality, it is enough to establish the upper bound of $\lambda^*_i$ for the proof of Lemma 3.5.

From (\ref{2.3}), (\ref{3.40}) and (\ref{3.39}), the flow (\ref{3.7}) has that for $u\in \mathbb{S}^{n-1}$ and $t\in (0, T]$,
\begin{equation}\label{3.43}
	\left\{
	\begin{aligned}
		\frac{\partial h^*}{\partial t}(u, t)&=-fe^{-\frac{1}{2h^*(u, t)^2}}\frac{\rho^*(x, t)^{n+1}}{h^*(u, t)^n\mathcal{S}^*}+(\sqrt{2\pi})^n\tau(t)h^*(u, t),\\
		h^*(u, t)&=\frac{1}{r},
	\end{aligned}
	\right.
\end{equation}
where 
\begin{eqnarray*}
\mathcal{S}^*=\det\left(\nabla^2h^* +h^*I \right)(u, t) 
\end{eqnarray*}
is the product of the principal radii of curvatures of $M^*_t$ at $p^*=\nabla h^*(u, t)+h^*(u, t)u$, and 
\begin{eqnarray*}
	\rho^*(x, t)=|p^*|=\sqrt{|\nabla h^*(u, t)|^2+h^*(u, t)^2}
\end{eqnarray*}
is the value of the radial function $\rho^*_{M_t^*}(\cdot, t)$ at 
\begin{eqnarray*}
x=\frac{p^*}{|p^*|}=\frac{\nabla h^*+h^*u}{\sqrt{|\nabla h^*|^2+(h^*)^2}}\in \mathbb{S}^{n-1}.
\end{eqnarray*}
Note also that $f$ takes its value at the above point $x$.

Let $b_{ij}^*=h_{ij}^*+h^*\delta_{ij}$, and let $b^{ij}_*$ be the inverse matrix of $b^*_{ij}$. As before, the eigenvalues of $b_{ij}^*$ and $b^{ij}_*$ are respectively the principal radii and principal curvatures of $M_t^*$. We consider the following auxiliary function
\begin{eqnarray*}
Q(u, t)=\log\lambda^*_{\max}\left( b_{ij}^*\right)-A\log h^*+B\left|\nabla h^* \right|^2 
\end{eqnarray*}
on $\mathbb{S}^{n-1}\times [0, T)$, where $\lambda^*_{\max}\left( b_{ij}^*\right)$ denotes the maximal eigenvalue of the matrix $b_{ij}^*$, and $A$ and $B$ are positive constants to be chosen later. We assume that $Q$ is attained its maximum at $(u_0, t_0)$, and take an orthonormal frame at the point such that $b_{ij}^*$ and $b_*^{ij}$ are diagonal and 
\begin{eqnarray*}
	\lambda_{\max}^*\left( b_{ij}^*(u_0, t_0)\right) =b_{11}^*(u_0, t_0).
\end{eqnarray*}
Then, we can write $Q$ as 
\begin{eqnarray}\label{3.44}
	Q(u, t)=\log b_{11}^*-A\log h^*+B\left|\nabla h^* \right|^2.
\end{eqnarray}
Hence, at this point $(u_0, t_0)$, we have
\begin{equation}\label{3.45}
	0=\nabla_iQ=b_*^{11}\nabla_ib^*_{11}-A\frac{\nabla_i h^*}{h^*}+2B\sum_k\nabla_kh^*\nabla_{ki}h^*,
\end{equation}
\begin{eqnarray}\label{3.46}
	0&\geq& \nabla_{ij}Q=b_*^{11}\nabla_{ij}b^*_{11}-\left( b_*^{11}\right)^2\nabla_i b_{11}^*\nabla_j b_{11}^*-A\left(\frac{\nabla_{ij}h^*}{h^*} -\frac{\nabla_{i}h^*\nabla_{j}h^*}{\left(h^*\right)^2}\right)  \nonumber \\
	&&+2B\sum_k \left(\nabla_{ki}h^*\nabla_{kj}h^*+\nabla_kh^*\nabla_{kij}h^* \right), 
\end{eqnarray}
\begin{eqnarray}\label{3.47}
	0&\leq &\frac{\partial Q}{\partial t}=b^{11}_*\partial_tb^*_{11}-A\frac{\partial_t h^*}{h^*}+2B\sum_k \nabla_kh^* \nabla_k(\partial_th^*) \nonumber\\ 
	&=&b^{11}_*\left(\nabla_{11} (\partial_t h^*)+\partial_t h^*\right)-A\frac{\partial_t h^*}{h^*}+2B\sum_k \nabla_kh^* \nabla_k(\partial_th^*). 
\end{eqnarray}
The flow (\ref{3.43}) can be written as
\begin{eqnarray}\label{3.48}
\log\left((\sqrt{2\pi})^n\tau h^*-\partial_t h^*\right)=-\log\mathcal{S}^*+\log\left( fe^{-\frac{1}{2(h^*)^2}}\frac{(\rho^*)^{n+1}}{(h^*)^n}\right). 
\end{eqnarray}
Let 
\begin{eqnarray*}
\Lambda(u, t)=\log\left( fe^{-\frac{1}{2(h^*)^2}}\frac{(\rho^*)^{n+1}}{(h^*)^n}\right). 
\end{eqnarray*}
Then differentiating (\ref{3.48}), at the point $(u_0, t_0)$ we have
\begin{eqnarray}\label{3.49}
\frac{(\sqrt{2\pi})^n\tau \nabla_kh^*-\nabla_k(\partial_t h^*)}{(\sqrt{2\pi})^n\tau h^*-\partial_t h^*}=-b_*^{ij}\nabla_kb_{ij}^*+\nabla_k\Lambda,
\end{eqnarray}
and
\begin{eqnarray}\label{3.50}
	\frac{(\sqrt{2\pi})^n\tau \nabla_{11}h^*-\nabla_{11}(\partial_t h^*)}{(\sqrt{2\pi})^n\tau h^*-\partial_t h^*}&=&\frac{\left( (\sqrt{2\pi})^n\tau \nabla_1 h^*-\nabla_1(\partial_t h^*)\right)^2 }{\left((\sqrt{2\pi})^n\tau h^*-\partial_t h^*\right)^2} \nonumber\\
&&-b_*^{ii}	\nabla_{11}b_{ii}^*+b^{ii}_*b_*^{jj}\left(\nabla_1 b_{ij}^* \right)^2+\nabla_{11}\Lambda. 
\end{eqnarray}
According to the Ricci identity, we see
\begin{equation*}
	\nabla_{11}b_{ii}^*=\nabla_{ii}b_{11}^*-b_{11}^*+b_{ii}^*.
\end{equation*}
Multiplying both sides of (\ref{3.50}) by $-b_*^{11}$, it follows that
\begin{eqnarray}\label{3.51}
	\frac{b_*^{11}\nabla_{11}(\partial_t h^*)-(\sqrt{2\pi})^n\tau b_*^{11}\nabla_{11}h^*}{(\sqrt{2\pi})^n\tau h^*-\partial_t h^*}&\leq&b_*^{11}b_*^{ii}	\nabla_{11}b_{ii}^*-b_*^{11}b^{ii}_*b_*^{jj}\left(\nabla_1 b_{ij}^* \right)^2-b_*^{11}\nabla_{11}\Lambda. \nonumber \\
	&\leq&b_*^{11}b_*^{ii}	\nabla_{ii}b_{11}^*-b_*^{11}b^{ii}_*b_*^{11}\left(\nabla_1 b_{i1}^* \right)^2-\sum_i b_*^{ii} \nonumber\\
	&&+b_*^{11}(n-1-\nabla_{11}\Lambda).
\end{eqnarray}
Note that $b_*^{ij}Q_{ij}\leq 0$ from (\ref{3.46}), which implies
\begin{eqnarray}\label{3.52}
&&b_*^{11}b_*^{ii}\nabla_{ii} b_{11}^*-(b_*^{11})^2b_*^{ii}(\nabla_ib_{11}^*)^2 \nonumber \\
 \ \ \ \ \ \ \ &\leq& Ab_*^{ii}\left(\frac{\nabla_{ii}h^*}{h^*} -\frac{(\nabla_{i}h^*)^2}{\left(h^*\right)^2}\right)-2B\sum_kb_*^{ii}\left((\nabla_{ki}h^*)^2+\nabla_{k}h^*\nabla_{kii}h^* \right).
\end{eqnarray}
We also need to calculate
\begin{eqnarray*}
b^{ii}_*\nabla_{ii}h^*=b_*^{ii}(b_{ii}^*-h^*)=n-1-h^*\sum_ib_*^{ii},
\end{eqnarray*}
\begin{eqnarray*}
\sum_kb_*^{ii}(\nabla_{ki}h^*)^2&=&b_*^{ii}(\nabla_{ii}h^*)^2=b_*^{ii}\left((b_{ii}^*)^2-2h^*b_{ii}^*+(h^*)^2 \right)\\
&=& -2(n-1)h^*+\sum_i b_{ii}^*+(h^*)^2\sum_i b_*^{ii},
\end{eqnarray*}
\begin{eqnarray*}
\sum_k b_*^{ii}\nabla_k h^* \nabla_{kii}h^*&=&\sum_kb_*^{ii}\nabla_k h^* \left( \nabla_i b_{ki}^*-\nabla_ih^*\delta_{ki}\right) \\
 &=&\sum_k b^{ii}_*\nabla_kh^*\nabla_{k}b_{ii}^*-b_*^{ii}(\nabla_i h^*)^2,
\end{eqnarray*}
where it is used that $\nabla_k b_{ij}^*$ is symmetric for all indices. Thus inequality (\ref{3.52}) can be calculated as
\begin{eqnarray*}
&&	b_*^{11}b_*^{ii}\nabla_{ii}b_{11}^*-(b_*^{11})^2b_*^{ii}(\nabla_i b_{11}^*)^2\leq A\frac{(n-1)}{h^*}-A\sum_i b_*^{ii}-A\frac{b_*^{ii}(\nabla_i h^*)^2}{(h^*)^2}-2B \sum_ib^*_{ii}\\
&&+4(n-1)Bh^*-2B(h^*)^2\sum_i b^{ii}_*-2B\sum_k b_*^{ii}\nabla_kh^*\nabla_kb_{ii}^*+2Bb_*^{ii}(\nabla_i h^*)^2.
\end{eqnarray*}
Plugging this inequality into (\ref{3.51}), we have
\begin{eqnarray}\label{3.53}
&&\frac{b_*^{11}\nabla_{11}(\partial_t h^*)-(\sqrt{2\pi})^n\tau b_*^{11}\nabla_{11}h^*}{(\sqrt{2\pi})^n\tau h^*-\partial_t h^*}	\leq A\frac{(n-1)}{h^*}-\left( A+2B(h^*)^2+1\right) \sum_ib_*^{ii}\nonumber \\
&&-\frac{A-2B(h^*)^2}{(h^*)^2}b_*^{ii}(\nabla_i h^*)^2+4(n-1)Bh^*-2B\sum_ib_{ii}^*-2B\sum_kb_*^{ii}\nabla_k h^*\nabla_k b_{ii}^*\nonumber \\
&& +b_*^{11}(n-1-\nabla_{11}\Lambda).
\end{eqnarray}
From (\ref{3.49}), we calculate
\begin{eqnarray}\label{3.54}
	\frac{2B\sum_k\nabla_kh^*\nabla_k(\partial_t h^*)}{(\sqrt{2\pi})^n\tau h^*-\partial_t h^*}&=&\frac{2(\sqrt{2\pi})^n\tau B|\nabla h^*|^2}{(\sqrt{2\pi})^n\tau h^*-\partial_t h^*}+2B\sum_kb_*^{ii}\nabla_kh^*\nabla_kb_{ii}^* \nonumber \\
&&-2B\langle\nabla h^*, \nabla \Lambda \rangle.
\end{eqnarray}
Dividing (\ref{3.47}) by $(\sqrt{2\pi})^n\tau h^*-\partial_t h^*$, we have
\begin{eqnarray*}
  0&\leq &\frac{b^{11}_*\left(\nabla_{11} (\partial_t h^*)+\partial_t h^*\right)}{(\sqrt{2\pi})^n\tau h^*-\partial_t h^*}-\frac{A\partial_t h^*}{h^*\left((\sqrt{2\pi})^n\tau h^*-\partial_t h^*\right)}+\frac{2B\sum_k \nabla_kh^* \nabla_k(\partial_th^*)}{(\sqrt{2\pi})^n\tau h^*-\partial_t h^*}\nonumber \\
  &=&\frac{b_*^{11}\left( \nabla_{11}(\partial_th^*)-(\sqrt{2\pi})^n\tau \nabla_{11}h^*+(\sqrt{2\pi})^n\tau b_{11}^*-(\sqrt{2\pi})^n\tau h^*+\partial_t h^*\right) }{(\sqrt{2\pi})^n\tau h^*-\partial_t h^*}\nonumber \\
\end{eqnarray*}
\begin{eqnarray}\label{3.55}
  &&-\frac{A\partial_t h^*}{h^*\left((\sqrt{2\pi})^n\tau h^*-\partial_t h^*\right)}+\frac{2B\sum_k \nabla_kh^* \nabla_k(\partial_th^*)}{(\sqrt{2\pi})^n\tau h^*-\partial_t h^*} \nonumber \\
	&=&\frac{b_*^{11}\left( \nabla_{11}(\partial_th^*)-(\sqrt{2\pi})^n\tau \nabla_{11}h^*\right) }{(\sqrt{2\pi})^n\tau h^*-\partial_t h^*}+\frac{2B\sum_k \nabla_kh^* \nabla_k(\partial_th^*)}{(\sqrt{2\pi})^n\tau h^*-\partial_t h^*}-b_*^{11}  \nonumber \\
&&	+\frac{A}{h^*}-\frac{(\sqrt{2\pi})^n\tau(A-1)}{(\sqrt{2\pi})^n\tau h^*-\partial_t h^*}.
\end{eqnarray}
Using (\ref{3.53}) and (\ref{3.54}), it follows from (\ref{3.55}) that
\begin{eqnarray}\label{3.56}
0&\leq& \frac{nA}{h^*}-\left(A+2B(h^*)^2+1\right)\sum_ib^{ii}_*-\frac{A-2B(h^*)^2}{(h^*)^2}\sum_ib_*^{ii}(\nabla_ih^*)^2 \nonumber \\
&&+4(n-1)Bh^*-2B\sum_ib_{ii}^*+b_*^{11}(n-2-\nabla_{11}\Lambda)-2B\langle\nabla h^*, \nabla \Lambda\rangle \nonumber \\
&&-(\sqrt{2\pi})^n\tau\frac{A-1-2B|\nabla h^*|^2}{(\sqrt{2\pi})^n\tau h^*-\partial_t h^*}
\end{eqnarray}
Recalling (\ref{3.41}), and noting (\ref{3.3}), we have
\begin{eqnarray}\label{3.57}
\left| \nabla h^*\right| \leq C,
\end{eqnarray}
for a positive constant $C$ independent of $t$. Thus taking $A=n+2BC^2$, we obtain
\begin{equation*}
A-1-2B\left| \nabla h^*\right|^2\geq A-1-2BC^2\geq 0.
\end{equation*}
Hence, it follows from (\ref{3.56}) that for some positive constant $C$,
\begin{eqnarray}\label{3.58}
 \ \ \ \ \ \ \ \  \ (A-n+3)\sum_ib_*^{ii}+2B\sum_ib_{ii}^*\leq C(A+B)-b^{11}_*\nabla_{11}\Lambda-2B\langle \nabla h^*, \nabla\Lambda\rangle.
\end{eqnarray}
Since
\begin{equation*}
\Lambda(u, t)=\log f-\frac{1}{2(h^*)^2}+(n+1)\log\rho^*-n\log h^*,
\end{equation*}
we calculate
\begin{eqnarray*}
	\nabla_k\Lambda&=&\frac{\nabla_k f}{f}-n\frac{\nabla_k h^*}{h^*}+\frac{\nabla_kh^*}{(h^*)^3}+(n+1)\frac{\nabla_k\rho^*}{\rho^*}\\
	&=&\frac{\nabla_k f}{f}-n\frac{\nabla_k h^*}{h^*}+\frac{\nabla_kh^*}{(h^*)^3}+(n+1)\frac{\nabla_k h^* b_{kk}^*}{(\rho^*)^2}.
\end{eqnarray*}
Thus,
\begin{eqnarray}\label{3.59}
	-2B\langle \nabla h^*, \nabla \Lambda\rangle&=&-2B\sum_k \left[\nabla_kh^*\left(\frac{\nabla_k f}{f}-n\frac{\nabla_k h^*}{h^*}+\frac{\nabla_kh^*}{(h^*)^3}+(n+1)\frac{\nabla_k h^* b_{kk}^*}{(\rho^*)^2}\right) \right]\nonumber \\
	&\leq& C B+C B\sum_k b_{kk}^*.
\end{eqnarray}
Moreover, 
\begin{eqnarray*}
-\nabla_{11}\Lambda& =&\frac{(\nabla_1 f)^2}{f^2}-\frac{\nabla_{11}f}{f}+n\frac{\nabla_{11}h^*}{h^*}-n\frac{(\nabla_1h^*)^2}{(h^*)^2}+\frac{\nabla_{11}h^*}{(h^*)^3}-\frac{3(\nabla_1h^*)^2}{(h^*)^4} \nonumber \\
&&-(n+1)\frac{\nabla_{11}h^*b_{11}^*}{(\rho^*)^2}-(n+1)\frac{\nabla_1h^*\nabla_1b_{11}^*}{(\rho^*)^2}+\frac{2(\nabla_1h^*)^2(b_{11}^*)^2}{(\rho^*)^4}\nonumber \\
&\leq& C\left( 1+b_{11}^*+(b_{11}^*)^2\right)+C\nabla_1b_{11}^*.
\end{eqnarray*}
Combining this with (\ref{3.45}), we obtain 
\begin{eqnarray}\label{3.60}
	-b^{11}_*\nabla_{11}\Lambda &\leq & C\left( 1+b^{11}_*+b_{11}^*\right)+Cb^{11}_*\nabla_1b_{11}^* \nonumber \\
	&=&C\left( 1+b^{11}_*+b_{11}^*\right)+C\left( A\frac{\nabla_1 h^*}{h^*}-2B\nabla_1h^*\nabla_{11}h^*\right). 
\end{eqnarray}
Adding both sides of inequalities (\ref{3.59}) and (\ref{3.60}), we have

\begin{eqnarray}\label{3.61}
-b^{11}_*\nabla_{11}\Lambda-2B\langle \nabla h^*, \nabla \Lambda\rangle &\leq& C\left( A\frac{\nabla_1 h^*}{h^*}-2B\nabla_1h^*\nabla_{11}h^*\right)+BC\sum_kb_{kk}^* \nonumber \\
&&+C\left( 1+b^{11}_*+b_{11}^*\right)+CB\nonumber \\
&\leq&Cb_*^{11}+Cb_{11}^*+C(A+B)+2BC\sum_ib_{ii}^* +C
\end{eqnarray}
Inserting (\ref{3.61}) into (\ref{3.58}), we have
\begin{eqnarray*}
(A-n+3)\sum_ib_*^{ii}+2B\sum_ib_{ii}^*\leq Cb_*^{11}+Cb_{11}^*+C(A+B)+2BC\sum_ib_{ii}^*+C,
\end{eqnarray*}
which has
\begin{eqnarray*}
	(A-n)\sum_ib_*^{ii}+2B\sum_ib_{ii}^*\leq C\sum_ib_*^{ii}+C\sum_ib_{ii}^*+C(A+B)+2BC\sum_ib_{ii}^*+C,
\end{eqnarray*}
Noting that $A=n+2BC^2$ as before, we get 
\begin{eqnarray*}
	(2BC^2-C)\sum_ib_*^{ii}+(2B-2BC-C)\sum_ib_{ii}^*\leq C(2BC^2+B+n)+C.
\end{eqnarray*}
If take
\begin{equation*}
B>\max\left\{\frac{C}{2C^2}, \frac{C}{2(1-C)}\right\}
\end{equation*}
in the above inequality, then $b^*_{11}(u_0, t_0)$ is bounded from above. From this, Lemma 3.3 and Corollary 3.1, and together with (\ref{3.44}), we easily see that the principal radii $\lambda_i^*$ ($i=1,\cdots, n-1$) of curvatures of $M_t^*$ have uniformly upper bound.\hfill $\square$

With the help of the above lemmas, we show that the flow (\ref{3.4}) exists for all time.

{\it \noindent{\bf Theorem 3.1.}~Let $f:\mathbb{S}^{n-1}\rightarrow (0, \infty)$ be a smooth, positive function. Then the flow (\ref{3.4}) has a smooth solution $M_t=X(\mathbb{S}^{n-1}, t)$ for all time $t>0$. Moreover, a subsequence of $M_t$ converges in $C^\infty$ to a smooth, closed, uniformly convex hypersurface, whose support function is a smooth solution to the equation:
	\begin{equation}\label{3.62}
fe^{-\frac{h^2+|\nabla h|^2}{2}} \mbox{\rm det}(\nabla_{ij}h+hI)=\frac{\int_{\mathbb{S}^{n-1}}e^{-\frac{h^2+|\nabla h|^2}{2}} (h^2+|\nabla h|^2)^{\frac{n}{2}}dx}{\int_{\mathbb{S}^{n-1}}h/fdx}.
\end{equation}}
\ \ {\bf Proof.}~By Lemma 3.3, Corollary 3.1 and Lemma 3.5, the equation (\ref{3.5}) is uniformly parabolic on any finite time interval. Thus, it follows from the standard parabolic theory and Krylov-Safonov estimates \cite{47a} that the smooth solution of the flow (\ref{3.5}) exists for all time.

Noting Lemma 3.3 and recalling definition (\ref{3.11}), there exists a positive constant $C$ which is independent of $t$, such that for any $t\geq 0$,
\begin{equation*}
\gamma_n(t)\leq C.
\end{equation*}
Thus,
\begin{equation*}
\int_0^t\gamma_n'(t)dt=\gamma_n(t)-\gamma_n(0)\leq\gamma_n(t)\leq C,
\end{equation*}
which leads to
\begin{equation*}
	\int_0^\infty\gamma_n'(t)dt\leq C.
\end{equation*}
According to Lemma 3.2, this inequality implies that there exists a subsequence of times $t_j\rightarrow \infty$ such that
\begin{equation*}
\gamma_n'(t_j)\rightarrow 0.
\end{equation*}
From this and Lemma 3.2, we see that there is a smooth funtion, say $h_\infty =h_{M_\infty}$, solving equation (\ref{3.62}). \hfill $\square$
\subsection{\bf Proof of Theorem 1.1}\hfill \break
\indent We use the approximation argument to prove Theorem 1.1.

{\bf Proof of Theorem 1.1.}~Let $f_j$ be a sequence of positive and smooth functions on $\mathbb{S}^{n-1}$ such that $\{\mu_j\}$, $d\mu_j=\frac{1}{f_j}dx$, converges weakly to $\mu$. For $f=f_j$, let $\{h_j\}$, $h_j=h_{K_j}$, be the solutions of the equation
	\begin{equation*}
	\frac{1}{(\sqrt{2\pi})^n}fe^{-\frac{h^2+|\nabla h|^2}{2}} \mbox{\rm det}(\nabla_{ij}h+hI)=\tau,
\end{equation*}
where
	\begin{equation*}
\tau=\frac{\int_{\mathbb{S}^{n-1}}e^{-\frac{h^2+|\nabla h|^2}{2}} (h^2+|\nabla h|^2)^{\frac{n}{2}}dx}{(\sqrt{2\pi})^n\int_{\mathbb{S}^{n-1}}h/fdx}.
\end{equation*}
Then 
\begin{equation*}
S_{\gamma_n, K_j}=\tau_j\mu_j~~\mbox{with}~~\tau_j=\frac{\int_{\mathbb{S}^{n-1}}e^{-\frac{h_j^2+|\nabla h_j|^2}{2}} (h_j^2+|\nabla h_j|^2)^{\frac{n}{2}}dx}{(\sqrt{2\pi})^n\int_{\mathbb{S}^{n-1}}h_j/f_jdx}.
\end{equation*}
Further, we show the lower and upper bounds of $h_{K_j}$ and $\tau_j$, which first need to prove that for a positive constant $c$,
\begin{equation}\label{3.63-}
	\min_{u\in \mathbb{S}^{n-1}}\left\| h_{\tilde{u}}\right\|_{\mu_j}\geq c>0.
\end{equation}
Note that $h_{\tilde{u}}\leq 1$ for any $u\in \mathbb{S}^{n-1}$. Then from (\ref{3.13}) we have $\left\| h_{\tilde{u}}\right\|_{\mu_j}\leq 1 $. By (\ref{3.20}), we also have that for sufficiently large $j$,
\begin{equation*}
	\min_{u\in \mathbb{S}^{n-1}}\left\|h_{\tilde{u}}\right\|_{\mu_j}>0.
\end{equation*}
Let 
\begin{equation}\label{3.63}
\lambda_j=\left\| h_{\tilde{u}_j}\right\| _{\mu_j}=	\min_{u\in \mathbb{S}^{n-1}}\left\|h_{\tilde{u}}\right\|_{\mu_j}
\end{equation}
for some $u_j\in \mathbb{S}^{n-1}$. We claim that $\{\lambda_j\}$ has the uniformly positive lower bound, and argue by contradiction. Suppose $\{\lambda_j\}\subset [0, 1]$
has a subsequence $\{\lambda_{j'}\}$ converging to $0$. From (\ref{3.63}), we can write that for some $u_{j'}\in \mathbb{S}^{n-1}$,
\begin{equation*}
\lambda_{j'}=\left\|h_{\tilde{u}_{j'}} \right\|_{\mu_{j'}}. 
\end{equation*}
Since $\mathbb{S}^{n-1}$ is compact, there exists a subsequence $\{\tilde{u}_{j^{''}}\}$ of $\{\tilde{u}_{j'}\}$ such that $\tilde{u}_{j^{''}}\rightarrow \tilde{u}_0\in \mathbb{S}^{n-1}$. Taking 
\begin{equation*}
	\lambda_{j^{''}}=\left\|h_{\tilde{u}_{j''}} \right\|_{\mu_{j''}},
\end{equation*}
we have $\lambda_{j^{''}}\rightarrow 0$ by the assumption. Note the fact that $h_{\tilde{u}}(x)=\frac{|\langle x, u\rangle|+\langle x, u\rangle}{2}$. Then $h_{\tilde{u}_{j^{''}}}\rightarrow h_{\tilde{u}_0}$ uniformly on $\mathbb{S}^{n-1}$. Thus from (\ref{2.1}) we get
\begin{equation*}
	\lambda_{j^{''}}=\left\|h_{\tilde{u}_{j''}} \right\|_{\mu_{j''}}\rightarrow \left\|h_{\tilde{u}_{0}} \right\|_{\mu}=0,
\end{equation*}
which contradicts to
\begin{equation*}
\left\|h_{\tilde{u}_{0}} \right\|_{\mu}\geq c>0
\end{equation*}
by (\ref{3.20}).

Therefore, from (\ref{3.19}) and (\ref{3.63-}) we see that $h_{K_j}$ has the uniformly upper bound. That is,
\begin{equation*}
	h_{K_j}\leq\frac{r}{\min_{x\in \mathbb{S}^{n-1}}\left\|h_{\tilde{u}}\right\|_{\mu_j}}\leq \frac{r}{c}.
\end{equation*}

Next, we derive the uniformly positive lower bound of $h_{K_j}$ by contradiction. It follows from Lemma 3.2 that
\begin{equation*}
	\gamma_n(K_t)=\gamma_n(t)\geq \gamma_n(0)=\frac{1}{2}.
\end{equation*}
Hence, if $K_\infty\in \mathcal{K}_o^n$ is a solution of the equation (\ref{3.62}) with $K_t\rightarrow K_\infty$, then from Lemma 3.3 and the dominated convergence theorem we have
\begin{equation*}
	\gamma_n(K_\infty)\geq\frac{1}{2}.
\end{equation*}
Thus by the assumption that $K_j$ are the solution of equation (\ref{3.62}), we obtain
\begin{equation}\label{3.65}
	\gamma_n(K_j)\geq\frac{1}{2}.
\end{equation}
We know that $K_j$ have the bound from above. Hence $K_i$ has a convergent subsequence, denoted also by $K_j$, which converges to a compact convex subset $K'\subset \mathbb{R}^n$. We will show that $K'$ has the origin in its interior. If not, then the origin is on the boundary of $K$. This implies that there exists $x'\in \mathbb{S}^{n-1}$ such that $h_{K'}(x')=0$. Since $K_i\rightarrow K'$ with respect to the Hausdorff metric, we get
\begin{equation*}
	\lim_{j\rightarrow \infty}h_{K_j}(x')=h_{K'}(x')=0.
\end{equation*}
Then for any $\varepsilon>0$ and sufficiently large $j$, we see $h_{K_j}(x')<\varepsilon$. That is,
\begin{equation}\label{3.66}
	K_j\subset \{y\in \mathbb{R}^n: \langle y, x'\rangle\leq \varepsilon\}.
\end{equation}
Noting that $K_j$ has upper bound, say $K_i\subset B_R$ for some $R>0$, and together with (\ref{3.66}) we have
\begin{equation*}
	K_j\subset B_R\cap\{y\in \mathbb{R}^n: \langle y, x'\rangle\leq \varepsilon\}.
\end{equation*}
It is well-known that
 \begin{eqnarray*}
	\gamma_n(\mathbb{R}^n)=\frac{1}{\left(\sqrt{2\pi}\right)^n}\int_{\mathbb{R}^n}e^{-\frac{|y|^2}{2}}dy=1.
\end{eqnarray*}
Then 
 \begin{eqnarray*}
	\gamma_n(H^-_{x'})=\frac{1}{\left(\sqrt{2\pi}\right)^n}\int_{H^-_{x'}\cap B_R}e^{-\frac{|y|^2}{2}}dy+\frac{1}{\left(\sqrt{2\pi}\right)^n}\int_{H^-_{x'}\setminus B_R}e^{-\frac{|y|^2}{2}}dy=\frac{1}{2}.
\end{eqnarray*}
Hence, for $\varepsilon$ small enough
\begin{eqnarray*}
	\gamma_n(K_j)\leq\gamma_n\left(B_R\cap\{y\in \mathbb{R}^n: \langle y, x'\rangle\leq \varepsilon\}\right)<\frac{1}{2},
\end{eqnarray*}
which is a contradiction to (\ref{3.65}). Namely, 
\begin{eqnarray*}
	h_{K_j}\geq c>0
\end{eqnarray*}
for a positive constant $c$.

By the Blaschke selection theorem, we obtain that $K_j$ subsequently converges to a convex body $K_0\in \mathcal{K}_o^n$. Thus, combining the weak convergence of the Gaussian surface area measure (see Section 2.1), it follows that
\begin{eqnarray*}
	S_{\gamma_n, K_0}=\tau \mu,
\end{eqnarray*}
where 
	\begin{equation*}
	\tau=\frac{\int_{\mathbb{S}^{n-1}}e^{-\frac{(\rho_{K_0})^2}{2}} (\rho_{K_0})^ndu}{(\sqrt{2\pi})^n\int_{\mathbb{S}^{n-1}}h_{K_0}d\mu}>0.
\end{equation*}
\hfill $\square$

\section{\bf Non-normalized solution}
\subsection{\bf Regularity of solution}\hfill \break
\indent This subsection is devoted to the study of the regularity of solutions to the Gaussian Minkowski problem, which is used to prove the existence of non-normalized and smooth solutions to the problem. Let us first recall some basic notions and facts required here. We refer to the papers \cite{48a} and \cite{49a} for some details.

  Let $K$ be a convex body. If $\partial K$ contains no segment then we say that $K$ is strictly convex; if $K$ has a unique tangential hyperplane at $p\in \partial K$ then we say that $p$ is a $C^1$-smooth point. Clearly, $h_K$ is $C^1$ on $\mathbb{S}^{n-1}$ if and only if $K$ is strictly convex. Moreover, $\partial K$ is $C^1$ if and only if each $p\in \partial K$ is $C^1$-smooth.
  
  Let $\Omega$ be a convex set in $\mathbb{R}^{n}$. If $y=\lambda z_1+(1-\lambda)z_2$ for $z_1, z_2\in \Omega$ and $\lambda\in (0, 1)$ implies $z_1=z_2=y$, then we say that $y\in \Omega$ is an extremal point. Note that if $\Omega$ is compact and convex, then $\Omega$ is the convex hull of its extremal points.
  
  The normal cone of a convex body $K$ at $z\in K$ is defined by
  \begin{eqnarray*}
  	N(K, z)=\left\lbrace x\in \mathbb{R}^{n}: \langle x, y\rangle\leq \langle x, z\rangle~\mbox{for all}~y\in K\right\rbrace, 
  \end{eqnarray*}
  or, equivalently by
  \begin{eqnarray*}
  	N(K, z)=\left\lbrace x\in \mathbb{R}^{n}:h_K(x)=\langle x, z \rangle\right\rbrace.
  \end{eqnarray*}
  If $z\in \mbox{int} K$ then we see $N(K, z)=\left\lbrace  o \right\rbrace$, and if $z\in \partial K$ then $\mbox{dim}N(K, z)\geq 1$.
  
  The face of a convex body $K$ with outer normal $x\in \mathbb{R}^{n}$ is defined as 
  \begin{eqnarray}\label{4.1}
  	F(K, x)=\left\lbrace z\in K:h_K(x)=\langle x, z \rangle\right\rbrace,
  \end{eqnarray}
  which lies in $\partial K$ provided $x\neq o$, and 
  \begin{eqnarray}\label{4.2}
  	F(K, x)=\partial h_K(x).
  \end{eqnarray}
  Here $\partial h_K(x)$ is the subgradient of $h_K$, namely
  \begin{eqnarray*}
  	\partial h_K(x)=\left\lbrace z\in \mathbb{R}^{n}: h_K(y)\geq h_K(x)+\langle z, y-x\rangle~\mbox{for each } y\in K\right\rbrace, 
  \end{eqnarray*}
  which is a non-empty compact convex set. Note that $h_K(x)$ is differentiable at $x$ if and only if $\partial h_K(x)$ consists of exactly one vector which is the gradient, denoted by $Dh_K(x)$, of $h_K$ at $x$.
  
 When $h_K$ is viewed as restricted to the unit sphere $\mathbb{S}^{n-1}$, the gradient of $h_K$ on $\mathbb{S}^{n-1}$ is denoted by $\nabla h_K$. Since $h_K$ is differentiable at $\mathcal{H}^n$ almost all points in $\mathbb{R}^n$ and is positively homogeneous of degree $1$, $h_K$ is differentiable for $\mathcal{H}^{n-1}$ almost all points of $\mathbb{S}^{n-1}$. Let $h_K$ be differentiable at $x\in \mathbb{S}^{n-1} $, where $x=\nu_K(p)$ is an outer unit normal vector at $p\in \partial K$. Then we have
\begin{eqnarray}\label{4.3}
	p=Dh_K(x)=\nu_K^{-1}(x).
\end{eqnarray}
This implies
\begin{eqnarray}\label{4.4}
	h_K(x)=h_K(\nu_K(p))=\langle p, \nu_K(p)\rangle=\langle Dh_K(x), x\rangle,
\end{eqnarray}
\begin{eqnarray}\label{4.5}
	Dh_K(x)=\nabla h_K(x)+h_K(x)x,
\end{eqnarray}
\begin{eqnarray}\label{4.6}
	|Dh_K(x)|=\sqrt{h_K(x)^2+|\nabla h_K(x)|^2}.
\end{eqnarray} 
 
The following is to recall the notions and facts about Monge-Amp\`{e}re measure from the survey by Trudinger and Wang \cite{50a}.

  Given a convex function $v$ defined in an open convex set $\Omega$ of $\mathbb{R}^{n}$, $Dv$ and $D^2v$ denote its gradient and its Hessian, respectively. For any Borel subset $\vartheta \subset\Omega$, define 
  \begin{eqnarray*}
  	N_v(\vartheta)=\bigcup_{y\in\vartheta}\partial v(y).
  \end{eqnarray*}
  The Monge-Amp\`{e}re measure $\mu_v$ is $\mu_v(\vartheta)=\mathcal{H}^n\left( N_v(\vartheta)\right)$. If the function $v$ is $C^2$ smooth, then the subgradient $\partial v$ coincide with the gradient $D v$. Thus
  \begin{eqnarray}\label{4.7}
  	\mu_v(\vartheta)=\mathcal{H}^n\left( Dv(\vartheta)\right)=\int_\vartheta \mbox{det}\left(D^2 v\right)d\mathcal{H}^n.
  \end{eqnarray}
  Note that the surface area measure $S_K$ of a convex body $K$ in $\mathbb{R}^{n}$ is a Monge-Amp\`{e}re type measure with $h_K$ restricted to the unit sphere $\mathbb{S}^{n-1}$ since it satisfies
  \begin{eqnarray}\label{4.8}
  	S_K(\omega)=\mathcal{H}^{n-1}\left(\bigcup_{x\in\omega}F(K, x)\right)=\mathcal{H}^{n-1}\left( \bigcup_{x\in\omega} \partial h_K(x)\right)=\mu_{h_K}(\omega)
  \end{eqnarray}
  for any Borel $\omega\subset \mathbb{S}^{n-1}$.
  
  We say that a convex function $v$ is the solution of a Monge-Amp\`{e}re equation in the sense of measure (or in the Aleksandrov sense), if it solves the corresponding integral formula for $\mu_v$ instead of the original formula for $\mbox{det}(D^2v)$.
  
  To study the regularity of the solution to the Gaussian Minkowski problem, we may first transfer the Monge-Amp\`{e}re type equation (\ref{1.3}) from $\mathbb{S}^{n-1}$ to $\mathbb{R}^{n-1}$. Therefore, we consider the restriction of a solution $h$ of (\ref{1.3}) to the hyperplane tangential to $\mathbb{S}^{n-1}$ at $e\in \mathbb{S}^{n-1}$.

  {\it \noindent{\bf Lemma 4.1.}~Let $e\in \mathbb{S}^{n-1}, K\in \mathcal{K}_o^n$, and $v: e^\perp\rightarrow \mathbb{R}$ with $v(y)=h_K(y+e)$. If $h=h_K$ is a solution of (\ref{1.3}) for non-negative functions $f$, then $v$ satisfies the standard Monge-Amp\`{e}re equation on $e^\perp$ in the sense of measure:
  	\begin{eqnarray}\label{4.9}
  		{\rm det} D^2v(y)=\frac{(\sqrt{2\pi})^n}{f(\pi(y))(\sqrt{1+|y|^2})^{n+1}}e^{\frac{\left| Dv(y)+(v(y)-\langle Dv(y), y\rangle)\cdot e\right|^2}{2}}.
  	\end{eqnarray}	
  }
  
  {\bf Proof.}~Assume that for $K\in \mathcal{K}_o^n$, $h=h_K$
  solves the equation (\ref{1.3}) and note that for $x\in \mathbb{S}^{n-1}$, 
   \begin{eqnarray}\label{4.10}
  dS(K, x)=\det\left(\nabla^2h_K(x)+h(x)I\right)d\mathcal{H}^{n-1}(x).
  \end{eqnarray}
Then it follows from (\ref{1.3})  and  (\ref{4.6}) that
  \begin{eqnarray}\label{4.11}
  	dS(K, x)=(\sqrt{2\pi})^ne^{\frac{|Dh_K(x)|^2}{2}}\frac{1}{f}d\mathcal{H}^{n-1}(x).
  \end{eqnarray}	
  
  For $e\in \mathbb{S}^{n-1}$, let $H_e$ be the hyperplane in $\mathbb{R}^{n}$ which is tangential to $\mathbb{S}^{n-1}$ at $e$, and $e^\perp$ the orthogonal complement of $\{re: r\in \mathbb{R}\}$ in $\mathbb{R}^n$. For $y\in e^\perp$, we have $y=\sum_{i=1}^{n-1}y_ie_i$, where $\{e_1, \cdots, e_{n-1}\}$ is a basis of $e^\perp$. Denoted by $\pi: e^\perp\rightarrow \mathbb{S}^{n-1}$ the radial projection from $H_e=e+e^\perp$ to $\mathbb{S}^{n-1}$, which is defined by
  \begin{eqnarray*}
  	\pi(y)=\frac{y+e}{\sqrt{1+|y|^2}}.
  \end{eqnarray*}	
  By the fact that 
  \begin{eqnarray}\label{4.12}
  	\langle \pi(y), e\rangle=\frac{1}{\sqrt{1+|y|^2}},
  \end{eqnarray}	
  we can calculate that the determinant of the Jacobin of the mapping $y\mapsto x=\pi (y)$ is
  \begin{eqnarray}\label{4.13}
  	\left| \mbox{Jac}\pi\right| =\left| \frac{dx}{dy}\right| =\left( \dfrac{1}{\sqrt{1+|y|^2}} \right)^n.
  \end{eqnarray}
  Let $v: e^\perp\rightarrow \mathbb{R}$ be the restriction of $h_K$ on $H_e$. That is,
  \begin{eqnarray}\label{4.14}
  	v(y)=h_K(y+e)=\sqrt{1+|y|^2}h_K(\pi(y)).
  \end{eqnarray}
  Then it is not hard to see by (\ref{4.1}) and (\ref{4.2}) that
  \begin{eqnarray}\label{4.15}
  	\partial v(y)= \partial h_K(y+e)\left| \right.e^\perp=F(K, y+e)\left|\right. e^\perp=F(K, \pi(y))\left| \right.e^\perp.
  \end{eqnarray}
  It follows from the homogeneity of degree $1$ and the differentiability of $h_K$ that
  \begin{eqnarray*}
  	Dh_K(y+e)=Dh_K(x),
  \end{eqnarray*}
  where $x=\pi (y)$. From this, we see
  \begin{eqnarray*}
  	Dv(y)=Dh_K(y+e)\left| \right.e^\perp=Dh_K(x)\left| \right.e^\perp.
  \end{eqnarray*}
  Therefore we can assume that
  \begin{eqnarray}\label{4.16}
  	Dh_K(x)=Dv (y)-q e
  \end{eqnarray}
  for some undetermined constant $q\in \mathbb{R}$. By (\ref{4.4}), we see
  \begin{eqnarray}\label{4.17}
  	h_K(x)=\langle Dh_K(x), x\rangle.
  \end{eqnarray}
  Note also that
  \begin{eqnarray}\label{4.18}
  	x=\pi(y)=\frac{y+e}{\sqrt{1+|y|^2}},~~h_K(x)=\dfrac{v(y)}{\sqrt{1+|y|^2}}.
  \end{eqnarray}
  Substituting (\ref{4.16}) and (\ref{4.18}) into (\ref{4.17}), it follows that
  \begin{eqnarray}\label{4.19}
  	q=\langle Dv(y), y\rangle-v(y).
  \end{eqnarray}
  Thus combining (\ref{4.16}) with (\ref{4.19}) we have
  \begin{eqnarray}\label{4.20}
  	Dh_K(x)=Dv(y)+\left(v(y)-\langle Dv(y), y\rangle\right) \cdot e.
  \end{eqnarray}
  For a Borel set $\vartheta \subset e^\perp$, based on (\ref{4.8}) it is easy to check that
  \begin{eqnarray}\label{4.21}
  	\mathcal{H}^{n-1}\left(\bigcup_{x\in \pi(\vartheta)}\left( F(K, x)\left| e^\perp\right. \right) \right)=\int_{\pi(\vartheta)}\langle x, e\rangle dS_K(x).
  \end{eqnarray}
  Thus it follows from (\ref{4.15}), (\ref{4.4}), (\ref{4.11}), (\ref{4.20}), (\ref{4.13}) and (\ref{4.14}) that
  \begin{eqnarray*}
  	&&\int_\vartheta \mbox{det}D^2 v (y)d\mathcal{H}^{n-1}(y)=\int_{\pi(\vartheta)}\langle x, e\rangle dS_K(x)\\
  	&=&(\sqrt{2\pi})^n\int_{\pi(\vartheta)}\langle x, e\rangle e^{\frac{|Dh_K(x)|^2}{2}}\frac{1}{f(x)}d\mathcal{H}^{n-1}(x)\\
  	&=&(\sqrt{2\pi})^n\int_\vartheta \frac{e^{\frac{\left| Dv(y)+(v(y)-\langle Dv(y), y\rangle)\cdot e\right|^2}{2}}}{f(\pi(y))(\sqrt{1+|y|^2})^{n+1}}d\mathcal{H}^{n-1}(y).
  \end{eqnarray*}
  This implies that $v$ satisfies (\ref{4.9}) on $e^\perp$.
 \hfill $\square$
  
  The next two lemmas by Caffarelli \cite{13a, 14a}, see also \cite{49a}, play crucial role in the proof of the regularity of solutions to the Gaussian Minkowski problem.
  
  {\it \noindent{\bf Lemma 4.2}(Caffarelli\cite{13a}).~Let $\lambda_2>\lambda_1>0$, and let $v$ be a convex function on an open bounded convex set $\Omega\subset \mathbb{R}^n$ such that
  	\begin{eqnarray*}
  		\lambda_1\leq \det D^2v \leq \lambda_2
  	\end{eqnarray*}	
  	in the sense of measure.
  	
  	\noindent(i) If $v$ is non-negative and $W=\{y\in \Omega: v(y)=0\}$ is not a point, then $W$ has no 
  	
  	extremal point in $\Omega$.
  	
  	\noindent(ii) If $v$ is strictly convex, then $v$ is $C^1$.
  }
  
  {\it \noindent{\bf Lemma 4.3}(Caffarelli\cite{14a}).~For real functions $v$ and $f$ on an open bounded convex set $\Omega\subset \mathbb{R}^n$, let $v$ be strictly convex, and let $f$ be positive and continuous such that 
  	\begin{eqnarray*}
  		\det D^2v =f
  	\end{eqnarray*}	
  	in the sense of measure.
  	
  	\noindent(i) Each $z\in \Omega$ has an open ball $B\subset \Omega$ around $z$ such that the restriction of $v$ to $B$ 
  	
  	is in $C^{1, \alpha}(B)$ for any $\alpha\in (0, 1)$.

  	\noindent(ii) If $f$ is in $C^\alpha (\Omega)$ for some $\alpha\in (0, 1)$, then each $z\in \Omega$ has an open ball $B\subset \Omega$ 
  	
  	around $z$ such that the restriction of $v$ to $B$ is in $C^{2, \alpha}(B)$.
  }
  
  With the help of the above lemmas, we prove the following results of regularity.
  
 {\it \noindent{\bf Theorem 4.1.}~Suppose that $d\mu=1/fd\mathcal{H}^{n-1}$ with $0<c_1\leq f \leq c_2$ on $\mathbb{S}^{n-1}$, and $K\in\mathcal{K}_o^n$ satisfies $dS_{\gamma_{n, K}}=\frac{1}{f}d\mathcal{H}^{n-1}$ on $\mathbb{S}^{n-1}$. Then
 	
 	\noindent(i) $\partial K$ is $C^1$ and strictly convex, and $h_K$ is $C^1$ on $\mathbb{R}^n\setminus\{o\}$;
 	
 	\noindent(ii) if $f$ is continuous, then the restriction of $h_K$ to $\mathbb{S}^{n-1}$ is in $C^{1, \alpha}$ for any $\alpha\in (0, 1)$;
 	
 	\noindent(iii) if $f \in C^\alpha(\mathbb{S}^{n-1})$ for $\alpha\in (0, 1)$, then $h_K$ is $C^{2, \alpha}$ on $\mathbb{S}^{n-1}$.}

  {\bf Proof.} For $e\in \mathbb{S}^{n-1}$ and $0<\tau< 1$, define $\Phi (e, \tau)$ by
  \begin{eqnarray*}
  	\Phi(e, \tau)=\{x\in\mathbb{S}^{n-1}: \langle x, e\rangle> \tau\}.
  \end{eqnarray*}
  Note that $h_K$ is continuous on $\mathbb{S}^{n-1}$ for $K\in \mathcal{K}_o^n$. Hence there exist $0<\tau_1<1$ and $\delta>0$ such that
  \begin{eqnarray*}
  	h_K(x)\geq\delta, ~\mbox{for}~x\in \mbox{cl}\Phi(e, \tau_1).
  \end{eqnarray*}
  Here, both $\tau_1$ and $\delta$ depend on $e$ and $K$. Moreover, there exists $0<\epsilon<1$ depending on $e$ and $K$ such that if some $x\in \mbox{cl}\Phi(e, \tau_1)$ is the outer normal at $p\in\partial K$, then
  \begin{eqnarray}\label{4.22}
  	\epsilon<|p|<\frac{1}{\epsilon}.
  \end{eqnarray}
  Recall that $\pi$ is defined by
  \begin{eqnarray*}
  	\pi(y)=\frac{y+e}{\sqrt{1+|y|^2}}, ~~y\in e^\bot.
  \end{eqnarray*}	
  Define 
  \begin{eqnarray*}
  	\Psi_e=\pi^{-1}(\Phi (e, \tau_1)).
  \end{eqnarray*}	
  Let $v: e^\perp\rightarrow \mathbb{R}$ satisfy the assumption of Lemma 4.1. Then for $y\in \Psi_e$, it follows that from (\ref{4.3}), (\ref{4.20}) and (\ref{4.22}) that
  \begin{eqnarray}\label{4.23}
  	\epsilon\leq |Dv(y)+\left( v(y)-\langle Dv(y), y\rangle\right) \cdot e|\leq \frac{1}{\epsilon}.
  \end{eqnarray}	
  Note that for $y\in \mbox{cl}\Psi_e$,
  \begin{eqnarray*}
  	v(y)=\sqrt{1+|y|^2}h_K\left(\dfrac{e+y}{\sqrt{1+|y|^2}} \right) \geq \delta.
  \end{eqnarray*}
  In addition, it is easy to see that $v$ also has an upper bound depending on $e$ and $K$ for $y\in \mbox{cl}\Psi_e$. Since it is assumed that for positive constants $c_1$ and $c_2$,  
  \begin{eqnarray*}
  	0<c_1\leq f \leq c_2.
  \end{eqnarray*}
  Thus it follows from Lemma 4.1 and (\ref{4.23}) that there exists $\zeta\in (0, 1)$ depending on $e$ and $K$ such that for $y\in \Psi_e$,
  \begin{eqnarray}\label{4.24}
  	\zeta\leq \det D^2v (y)\leq \frac{1}{\zeta}.
  \end{eqnarray}
  
  We first show that $\partial K$ is $C^1$ for $K\in \mathcal{K}_o^n$. Namely, $\mbox{dim} N(K, z)=1$ for any 
  $z\in \partial K$. If not, then there exists a point $z_0\in \partial K$ such that $\mbox{dim} N(K, z_0)\geq 2$. Assume $e\in N(K, z_0) \cap \mathbb{S}^{n-1}$. By the definition of support function, and noting $z_0\in \partial K$, it follows that for $y\in \Psi_e$,
  \begin{eqnarray*}
  	v(y)=h_K(y+e)\geq \langle y+e, z_0\rangle.
  \end{eqnarray*}
  It is not hard to verify that 
  \begin{eqnarray*}
  	y\in W:=\pi^{-1}(N(K, z_0)\cap \Phi (e, \tau_1))\Longleftrightarrow v(y)=\langle y+e, z_0\rangle~\mbox{for}~y \in \Psi_e.
  \end{eqnarray*}
  Let $\varphi(y)=\langle y+e, z_0\rangle$. Then
  \begin{eqnarray*}
  	v(y)-\varphi(y)\left\{
  	\begin{aligned}
  		&=&0~~\mbox{for}~&y\in W&\\
  		&>&0~~\mbox{for}~&y\in \Psi_e\setminus W.&
  	\end{aligned}
  	\right.	
  \end{eqnarray*}
  By (\ref{4.24}), and noting that $\varphi$ is the first degree polynomial, we get that for $y\in \Psi_e$,
  \begin{eqnarray*}
  	\zeta\leq \det D^2(v (y)-\varphi(y))\leq \frac{1}{\zeta}.
  \end{eqnarray*}
  Note that
  \begin{eqnarray*}
  	W=\pi^{-1}(N(K, z_0)\cap \Phi (e, \tau_1))=\{y\in \Psi_e: v (y)-\varphi(y)=0\}
  \end{eqnarray*}
  is not a point since $\mbox{dim} W\geq 1$. Since the origin $o$ is an extremal point of $W$ by the choice of $e$, this is a contradiction from (i) of Lemma 4.2. 
  
  Next, we prove that $v$ is strictly convex on $\mbox{cl} \Psi_e$ for $e\in \mathbb{S}^{n-1}$. It is easy to see that $\Psi_e$ is a convex set in $\mathbb{R}^{n-1}$. For $0<\lambda<1$ and $y_1, y_2\in \Psi_e$ with $y_1\neq y_2$, we assume that $e+(\lambda y_1+(1-\lambda)y_2)$ is an outer normal at $z\in \partial K$. Namely,
  \begin{eqnarray*}
  	e+(\lambda y_1+(1-\lambda)y_2)\in N(K, z).
  \end{eqnarray*}
  Note that
  \begin{eqnarray*}
  	e+y_1\notin N(K, z)~~\mbox{and}~~e+y_2\notin N(K, z)
  \end{eqnarray*}
  since $z\in \partial K$ is a smooth point. Thus,
  \begin{eqnarray*}
  	v(y_i)=h_K(e+y_i)>\langle z,e+y_i\rangle.
  \end{eqnarray*}
  From this, we deduce
  \begin{eqnarray*}
  	\lambda v(y_1)+(1-\lambda) v (y_2)&>&
  	\langle z, e+(\lambda y_1+(1-\lambda)y_2)\rangle\\
  	&=&h_K(e+(\lambda y_1+(1-\lambda)y_2))\\
  	&=&v(\lambda y_1+(1-\lambda)y_2).
  \end{eqnarray*}
  That is,
  \begin{eqnarray*}
  	v(\lambda y_1+(1-\lambda)y_2)<\lambda v(y_1)+(1-\lambda) v (y_2).
  \end{eqnarray*}
  Therefore, according to (ii) of Lemma 4.2, it follows from (\ref{4.24}) and the strict convexity of $v$ that $v$ is $C^1$ on $\Psi_e$ for any $e\in \mathbb{S}^{n-1}$. This implies that $h_K$ is $C^1$ on $\mathbb{R}^{n}\setminus\{o\}$, and $\partial K$ contains no segment. This completes the proof of (i) in Theorem 4.1.
  
  We next prove (ii) in Theorem 4.1. From the assumption that $f$ is continuous, and noting that $v$ is $C^1$ on $\mbox{cl} \Psi_e$ for any $e\in \mathbb{S}^{n-1}$, it follows that the right hand side of (\ref{4.9}) is continuous. Using (i) of Lemma 4.3 and the strict convexity of $v$ on $\Psi_e$, we have that there exists an open ball $B\subset \Psi_e$ centred at $o$ such that $v$ is $C^{1, \alpha}$ on $B$ for any $\alpha\in (0, 1)$. This implies that $h_K$ is locally $C^{1, \alpha}$ on $\mathbb{S}^{n-1}$. Thus it follows from the compactness of $\mathbb{S}^{n-1}$ that $h_K$ is globally $C^{1, \alpha}$ on $\mathbb{S}^{n-1}$, which finishes the proof of (ii) in Theorem 4.1.
  
  Now let us show (iii) of Theorem 4.1. Since $v$ is $C^{1, \alpha}$ on $B$, together with the assumption that $f$ is $C^\alpha$ on $\mathbb{S}^{n-1}$ we obtain that the right hand side of (\ref{4.9}) is $C^\alpha$. According to (ii) of Lemma 4.3, we deduce that $v$ is $C^{2, \alpha}$ on an open ball $\widetilde{B}\subset B$ of $e^\perp$ centred at $o$. This implies that $h_K$ is locally $C^{2, \alpha}$ on $\mathbb{S}^{n-1}$. Thus $h_K$ is globally $C^{2, \alpha}$ on $\mathbb{S}^{n-1}$ based on the compactness of $\mathbb{S}^{n-1}$. This gives the proof of (iii) in Theorem 4.1.\hfill $\square$
  
 \subsection{\bf Existence of smooth solution}\hfill \break
 \indent Here, we use the degree theory to  obtain the existence of smooth solutions to the Gaussian Minkowski problem. The proof is based on the idea of Huang, Xi and Zhao \cite{1a}. Concerning the degree theory for second-order nonlinear elliptic operators, the reader may wish to consult the work of Li \cite{51a}.

The following lemma established in \cite{12a} will be needed, which contains a uniform estimate.

 {\it \noindent{\bf Lemma 4.4 \cite{12a}.}~Let $1/C<f<C$ for some $C>0$ and $K\in \mathcal{K}_o^n$ with $\gamma_n(K)>\frac{1}{2}$.  If $h_K\in C^2(\mathbb{S}^{n-1})$ is a solution of (\ref{1.3}), then there exists a constant $C'>0$ such that $1/C'<h_K<C'$.
 }

The existence result of smooth solution to the Gaussian Minkowski problem may be formulated as follows:
 
 {\it \noindent{\bf Theorem 4.2.}~Let $f$ be a positive smooth function on $\mathbb{S}^{n-1}$ and $|\frac{1}{f}|_{L_1}<\frac{1}{\sqrt{2\pi}}$. Then there exists a smooth convex body $K\in \mathcal{K}_o^n$ with $\gamma_n(K)>\frac{1}{2}$ such that $h_K$ satisfies (\ref{1.3}).
}

{\bf Proof.}~Let $g=\frac{1}{f}\equiv c_0>0$ small enough and $|c_0|_{L_1}<\frac{1}{\sqrt{2\pi}}$. Then it follows from the intermediate value theorem and Lemma 2.1 that there exists a unique constant solution $h_K=s_0$ of equation (\ref{1.3}) such that $\gamma_n(K)>\frac{1}{2}$. We consider a family of equations
\begin{eqnarray}\label{4.25}
	\frac{1}{(\sqrt{2\pi})^n}e^{-\frac{|\nabla h|^2+h^2 }{2}}\det\left(\nabla^2h+hI\right)=g_t ~~\mbox{on}~~ \mathbb{S}^{n-1},
\end{eqnarray}
where $t\in [0, 1]$, and $g_t=(1-t)c_0+tg$. We use $C_c^{2, \alpha}(\mathbb{S}^{n-1})$ to denote a set of positive and strictly convex functions from $C^{2, \alpha}(\mathbb{S}^{n-1})$ with $\alpha\in (0, 1)$. Define a family of operators $\mathcal{G}_t: C_c^{2, \alpha}(\mathbb{S}^{n-1})\rightarrow C^{\alpha}(\mathbb{S}^{n-1})$ by
\begin{eqnarray*}
\mathcal{G}_t(h)=}\det\left(\nabla^2h+hI\right)-(\sqrt{2\pi})^ne^{\frac{|\nabla h|^2+h^2 }{2}g_t.
\end{eqnarray*}
Let $g\in C^{\alpha}(\mathbb{S}^{n-1})$ and $g>0$. Then there exists a positive constant $C$ such that $\frac{1}{C}<g, c_0< C$ and $|g|_{C^{\alpha}(\mathbb{S}^{n-1})}<C$. Thus, $g_t$ also satisfy $\frac{1}{C}<g_t< C$, $|g_t|<\frac{1}{\sqrt{2\pi}}$ and $|g_t|_{C^{\alpha}}<C$. From Theorem 4.1 and Lemma 4.4, we can construct an open bounded subset $\mathcal{O}$ of $C_c^{2, \alpha}(\mathbb{S}^{n-1})$,
\begin{eqnarray*}
	\mathcal{O}=\left\lbrace h\in C_c^{2, \alpha}(\mathbb{S}^{n-1}): \frac{1}{C_1}<h<C_1, |h|_{C^{2, \alpha}}<C_2, \gamma_n([h])>\frac{1}{2}\right\rbrace, 
\end{eqnarray*}
where the positive constant $C_1$ depends on $C$ in Lemma 4.4, and $C_2$ can be obtained from Theorem 4.1. 

Next, we show by contradiction that if $h\in \partial \mathcal{O}$, then $\mathcal{G}_t(h)\neq 0$. Assume $\mathcal{G}_t(h)=0$. Namely, $h$ solves the equation (\ref{4.25}). Thus by Lemma 4.4 and Theorem 4.1 we only have $\gamma_n([h])=\frac{1}{2}$. Hence, from Lemma 2.2 it follows that $|S_{\gamma_n, K}|\geq \frac{1}{\sqrt{2\pi}}$, which contradicts with the fact that if $h$ solves the equation (\ref{4.25}) with $|g_t|_{L_1}<\frac{1}{\sqrt{2\pi}}$ then $|S_{\gamma_n, K}|< \frac{1}{\sqrt{2\pi}}$. Therefore,
 \begin{equation*}
	\mathcal{G}_t^{-1}(0)\cap \partial\mathcal{O}=\emptyset,
\end{equation*}
for all $t\in [0, 1]$. That is, the degree $\deg (	\mathcal{G}_t, \mathcal{O}, 0)$ is well-defined for all $0\leq t\leq 1$ and is independent of $t$. If we verify $\deg(\mathcal{G}_1, \mathcal{O}, 0)\neq 0$, then there exists a $h\in \mathcal{O}$ such that $\mathcal{G}_1(h)=0$, which implies that $h$ solves the equation (\ref{1.3}). Hence, from the fact that $\deg (	\mathcal{G}_t, \mathcal{O}, 0)$ is independent of $t$, i.e.,
 \begin{equation*}
	\deg(\mathcal{G}_1, \mathcal{O}, 0)=\deg(\mathcal{G}_0, \mathcal{O}, 0),
\end{equation*}
we just need to prove
\begin{equation*}
\deg(\mathcal{G}_0, \mathcal{O}, 0)\neq 0.
\end{equation*}

Let $\mathcal{L}_{s_0}: C^{2, \alpha}(\mathbb{S}^{n-1})\rightarrow C^{\alpha}(\mathbb{S}^{n-1})$
denote the linearized operator of $\mathcal{G}_0$ at $s_0$. It is not hard to compute
\begin{equation*}
	\mathcal{L}_{s_0}(\phi)=s_0^{n-2}\left(\triangle_{\mathbb{S}^{n-1}}\phi+((n-1)-s_0^2)\phi\right).
\end{equation*}
Since spherical Laplacian is a discrete spectrum, we can find a small enough $c_0$ with $|c_0|_{L_1}<\frac{1}{\sqrt{2\pi}}$ such that the operator $\mathcal{L}_{s_0}$ is invertible. According to the Propositions 2.3 and 2.4 in Li \cite{51a}, and combining with the fact that $s_0$ is the unique solution of $\mathcal{G}_0(h)=0$, we have
 \begin{equation*}
	\deg(\mathcal{G}_0, \mathcal{O}, 0)=\deg(\mathcal{L}_{s_0}, \mathcal{O}, 0)\neq 0.
\end{equation*}

By differentiating equation (\ref{4.25}) with $t=1$ repeatedly, we obtain the existence of smooth solution to equation (\ref{1.3}).  \hfill $\square$

 \subsection{\bf Proof of Theorem 1.2}\hfill \break
\indent This subsection is devoted to the proof of Theorem 1.2 by using an approximation argument.

{\bf Proof of Theorem 1.2.}~Let $d\mu_j=\frac{1}{f_j}dx$ with $0<f_j\in C^{\infty}(\mathbb{S}^{n-1})$ and $0<|\frac{1}{f_j}|_{L_1}<\frac{1}{\sqrt{2\pi}}$, and let $\{\mu_j\}$ converge weakly to $\mu$ as $j\rightarrow \infty$. It follows from Theorem 4.2 that there is a sequence of convex bodies $K_j$ with $\gamma_n(K_j)>\frac{1}{2}$ such that $h_j=h_{K_j}$ satisfy the equation
 \begin{eqnarray}\label{4.26-}
	\frac{1}{(\sqrt{2\pi})^n}e^{-\frac{|\nabla h_{K_j}|^2+h^2_{K_j} }{2}}\det(\nabla^2h_{K_j}+h_{K_j}I)=\frac{1}{f_j}.
\end{eqnarray}
Thus,
\begin{eqnarray*}
S_{\gamma_n, K_j}=\mu_j.
\end{eqnarray*}
Let $h_{K_j}(x_j)=\max_{x\in\mathbb{S}^{n-1}}h_{K_j}(x)$ for some $x_j\in\mathbb{S}^{n-1}$. Since the support of $\mu$ is not contained in a closed hemisphere, it follows by Lemma 5.4 in \cite{12a} that there exist positive constant $c$ and $\xi$ such that for large $j$,
\begin{eqnarray}\label{4.26}
\int_{\{x\in \mathbb{S}^{n-1}:\langle x, x_j\rangle>\xi\}}\frac{1}{f_j}dx\geq c>0.
\end{eqnarray}
From (\ref{2.5+}), we have that for any $x\in \mathbb{S}^{n-1}$ with $\langle x, x_j\rangle\geq \xi$,
\begin{eqnarray*}
h_{K_j}(x)\geq\langle x, x_j\rangle h_{K_j}(x_j)\geq \xi h_{K_j}(x_j).
\end{eqnarray*}
Note also that 
\begin{eqnarray*}
\det(\nabla^2h_{K_j}+h_{K_j}I)dx=\frac{(\rho_{K_j})^n}{h_{K_j}}du,
\end{eqnarray*}
where $x$ and $u$ are associated by
\begin{eqnarray*}
\rho_{K_j}(u)u=h_{K_j}(x)x+\nabla h_{K_j}(x).
\end{eqnarray*}
 Thus, it follows from (\ref{2.4}), (\ref{4.26-}) and (\ref{4.26}) that
\begin{eqnarray*}
	0<c&\leq &\int_{\{x\in \mathbb{S}^{n-1}:\langle x, x_j\rangle>\xi\}}\frac{1}{f_j}dx \\
	&\leq &	\frac{1}{(\sqrt{2\pi})^n}e^{-\frac{(\xi h_{K_j}(x_j))^2 }{2}}(h_{K_j}(x_j))^{n-1}\frac{1}{\xi}\mathcal{H}^{n-1}(\mathbb{S}^{n-1})\rightarrow 0,
\end{eqnarray*}
as $h_{K_j}(x_j)\rightarrow \infty$. This is a contradiction. Therefore, $K_j$ is uniformly bounded from above.

The lower bound of $K_j$ is guaranteed by $\gamma_n(K_i)>\frac{1}{2}$ and we see Section 3.3 for its proof. 

In conclusion, there exists a positive constant $C$ such that 
\begin{eqnarray}\label{4.28}
\frac{1}{C}\leq h_{K_j}\leq C.
\end{eqnarray}
By (\ref{4.28}) and the Blaschke selection theorem, there is a convergent subsequence $K_{j_i}$ of $K_j$ such that 
\begin{eqnarray*}
	K_{j_i}\rightarrow K_0\in \mathcal{K}_o^n
\end{eqnarray*}
in the Hausdorff metric and $\gamma_n(K_0)>\frac{1}{2}$. Hence, together with the weak convergence of the Gaussian surface area measure (see Section 2.1), we have
\begin{eqnarray*}
	S_{\gamma_n, K_0}=\mu,
\end{eqnarray*}
which completes the proof of Theorem 1.2. \hfill $\square$

\vskip10pt

\bibliographystyle{amsplain}

\end{document}